\documentclass[final, 3p, times,review]{elsarticle}

\usepackage{lineno,hyperref}

\usepackage{amsthm}
\usepackage{amsmath,amssymb,amsfonts,amscd}
\usepackage{latexsym}
\usepackage{algorithmicx,algorithm}
\usepackage{subfigure}
\usepackage{epstopdf}

\usepackage{url}
\usepackage{xcolor}
\definecolor{newcolor}{rgb}{.8,.349,.1}

\newcommand{\beq}{\begin{equation}}
\newcommand{\eeq}{\end{equation}}
\newcommand{\bal}{\begin{align}}
\newcommand{\eal}{\end{align}}
\newcommand{\baln}{\begin{align*}}
\newcommand{\ealn}{\end{align*}}

\theoremstyle{plain}
\theoremstyle{definition}\theoremstyle{plain}
\newtheorem{pro}{Proposition}

\graphicspath{{figure/}}

\begin{document}

\begin{frontmatter}

\title{Hamiltonian Particle-in-Cell methods for Vlasov--Poisson equations\tnoteref{tnote1}}%

\author[1,2]{Anjiao Gu}
\cortext[cor1]{Corresponding author at:
LSEC, ICMSEC, Academy of Mathematics and Systems Science, Chinese Academy of Sciences, Beijing 100190, China; School of Mathematical Sciences, University of Chinese Academy of Sciences, Beijing 100049, China.}
\author[3]{Yang He}
\author[1,2]{Yajuan Sun\corref{cor1}}
\ead{sunyj@lsec.cc.ac.cn}

\address[1]{LSEC, ICMSEC, Academy of Mathematics and Systems Science, Chinese Academy of Sciences, Beijing 100190, China}
\address[2]{School of Mathematical Sciences, University of Chinese Academy of Sciences, Beijing 100049, China}
\address[3]{School of Mathematics and Physics, University of Science and Technology Beijing, Beijing 100083, China}

\begin{abstract}
In this paper, Particle-in-Cell algorithms for the Vlasov--Poisson system are presented based on its Poisson bracket structure. The Poisson equation is solved by finite element methods, in which the appropriate finite element spaces are taken  to guarantee that  the semi-discretized system possesses a well defined  discrete Poisson bracket structure.
Then, splitting methods are applied to the semi-discretized system by decomposing  the Hamiltonian function. The resulting discretizations are proved to be Poisson bracket preserving. Moreover, the conservative quantities of the system are also well preserved.
In numerical experiments, we use the presented  numerical methods to  simulate various physical phenomena. Due to the huge computational effort of the practical computations, we employ the strategy of parallel computing. The numerical results verify the efficiency of the new derived numerical discretizations.
\end{abstract}


\end{frontmatter}

\section{Introduction}

The motion of charged particles under electromagnetic fields is the most fundamental physical process in magnetized plasma. If only one particle is concerned, the dynamics can be described by the single particle model. If we are researching the interactive dynamics of a large number of particles with self-consistent electromagnetic fields, the kinetic model can be applied. In this model the Vlasov equation is applied together with the Maxwell equations or the Poisson equation, and the coupled  system is called the Vlasov--Maxwell (VM) system or Vlasov--Poisson (VP) system respectively.

There are two main classes of numerical methods to solve the Vlasov type equations, the Eulerian method and the particle method. The Eulerian method, also called the grid-based method, is to discretize  the PDEs model on fixed computational grid, and to carry out the time integration of the distribution function on mesh points. Various Eulerian approaches such as finite element methods, finite volume methods and  spectral methods etc. have been developed and applied to the  plasma problems, see~\citep{Filbet2003N,Morrison2016S,XJY2018S} and references therein.
Different from the grid-based method the Particle-in-Cell (PIC) approach is to approximate the distribution function by `super particles' with the weighted Klimontovich representation, and to follow the trajectories of super particles.
Indeed, PIC has a relatively small computational cost to conduct high dimensional simulations~\citep{Nicolas2009,LYZ2019}.
Although the introduction of super particles brings computational noise, it has been confirmed that the computational noise will decrease in the rate of $1/\sqrt{N_p}$ along the increasing number $N_p$ of particles.
The numerical simulation using the PIC technique can reproduce well the  realistic physical phenomena~\citep{Hockney1988,Birdsall1991}, and has been considered as an effective way to simulate plasmas of the kinetic theory~\citep{Cottet2000}.

In most applications, there span a wide range of space and time scales for the problem, which demands numerical simulations over large time intervals. Classcial numerical methods, such as RK4, can approximate the solution well in the first few number of iterations.
However, the numerical error accumulates rapidly over steps, and usually leads to wrong numerical results after a large number of iterations~\citep{FK1995}. One reason is that  these methods can not preserve the important conservative  quantities of the original system.
Structure-preserving algorithms (or geometric numerical integrators) are designed to
conserve the intrinsic properties of the original system, including the symplectic structure, the Poisson structure, the invariant phase space volume, and constants of motions etc.. As a result the long-term stability of the numerical results can be guaranteed~\citep{FK1985,Hairer2006,FK2010}.
In~\citep{He2017Explicit,He2016High}, the $K$-symplectic algorithms and the volume-preserving algorithms have been developed for solving the single particle equations. These numerical algorithms provide the numerical simulations  with bounded energy error, and the trajectories of particles are well simulated over exponentially long simulation time.

It is known that the Vlasov--Maxwell system can be written  as a  Hamiltonian system with respect to the Morrison-Marsden-Weinstein (MMW) Poisson bracket~\citep{Morrison1980T,Marsden1982T}. The Vlasov--Poisson system can be taken as the reduced Vlasov--Maxwell system. In this paper, we study the Poisson bracket structure of the Vlasov--Poisson system, and verify that the Poisson bracket of the Vlasov--Poisson system has very closed relation to the MMW Poisson bracket. As the Poisson bracket is usually not canonical, conventional time integrators are not structure-preserving. However, the splitting technique~\citep{Yoshida1990,McLachlan2002} can be applied. High-order Hamiltonian splitting methods have been developed for Vlasov--Maxwell equations ~\citep{XJY2015E,Qin2016C,HY2016,Kraus2017,LYZ2019} and Vlasov--Poisson equations without magnetic fields~\citep{Nicolas2015H,Casas2017}.
Correspondingly, spatial discretizations should also be designed carefully to maintain the structure.  For the VM system, with the help of specially designed finite element spaces, the researchers~\citep{HY2016,Kraus2017,Burby2017Finite} provide the semi-discrete system which can conserve the discrete Poisson structure  corresponding to the original one. In other ways, a discrete Poisson bracket can also be obtained by applying discrete exterior calculus~\citep{XJY2015E,Qin2016C,Kraus2017} and the variational principle~\citep{Squire2012G,Evstatiev2013V,Martin2021Variational}.
On the other hand, in simulations of magnetic confinement fusion, Vlasov equation with a strong non-homogeneous magnetic field has aroused wide interest recently. Asymptotic preserving techniques~\citep{JS1999E,JS2017T,Filbet2016,Hairer2021L} and uniformly accurate schemes~\citep{BWZ2014U,Chartier2020} are efficient and characteristic numerical methods for them.

In this paper, we have developed the Hamiltonian Particle-in-Cell methods for the VP system with an external magnetic field. We first derive particle equations from the discretisation of the Vlasov equation, then the Poisson equation is solved by finite element methods, with the finite element spaces chosen such that the semi-discrete system is Hamiltonian with a discrete Poisson structure. Furthermore, the Hamiltonian of the semi-discrete system is split into two parts, each of which possesses the same Poisson structure and can be solved exactly. By composing the  solutions to the subsystems results in the  Poisson-structure-preserving discretizations for the semi-discrete equations. In the numerical  simulations, the parallel technique has been employed, and the effectiveness of the resulting numerical discretizations has been tested
in numerical experiments.

The outline of the paper is as follows. In section 2, we introduce the various descriptions for the  Vlasov equation. For the VP equation, in section 3 we  present the Hamiltonian formulation which is described with the defined  Poisson bracket structure. Furthermore, we apply the finite element methods to construct the numerical methods  for VP system in Section 4. The corresponding discrete Poisson bracket and the discrete
Hamiltonian are also established  in this section. In Section 5, we split the Hamiltonian in order to construct the numerical discretizations which can preserve the dicrete Poisson bracket. We display the numerical results in Section 6. Finally, we conclude this paper.

\section{Equations}
In the kinetic model, the Vlasov equation is in the form
$$\frac{\partial f}{\partial t}+\mathbf{v}\cdot\frac{\partial f}{\partial\mathbf{x}}+\frac{\mathbf{F}(\mathbf{x},\mathbf{v},t)}{m}\cdot\frac{\partial f}{\partial\mathbf{v}}=0,$$
where $f(\mathbf{x},\mathbf{v},t)$ is the distribution function of position $\mathbf{x}\in\Omega_{x}\subset\mathbb{R}^{3}$
and velocity $\mathbf{v}\in\Omega_{v}\subset\mathbb{R}^{3}$ at time $t$, and $\mathbf{F}$ is the force field acting on the particles.
In most cases, the charged particles are considered to interact with the self-consistent electromagnetic fields in vacuum. This yields the following Maxwell's equations,
\begin{equation*}
\begin{gathered}
\mu_0\epsilon_0\frac{\partial\mathbf{ E}}{\partial t}=\nabla\times\mathbf{ B}-\mu_0q\int_{\Omega_v}\mathbf{v} fd\mathbf{v}, \quad \nabla\cdot\mathbf{B}=0,\\
\frac{\partial\mathbf{B}}{\partial t}=-\nabla\times\mathbf{E}, \quad \epsilon_0\nabla\cdot\mathbf{E}=q\int_{\Omega_v}fd\mathbf{v}-q\rho_0,
\end{gathered}
\end{equation*}
where $\mu_0$ is the magnetic permeability, $\epsilon_0$ is the vacuum permittivity and the constant $\rho_0$ is the charge density of ions. Then the force is given by $\mathbf{F}=q(\mathbf{E+v\times B})$ with $q$ the charge.
\begin{table}[h]
\centering
\begin{tabular}{|l|c|c|}
  \hline
  Names & Symbols & ~~~Units~~~\\\hline
  Time & $t$ & $1/\omega_p$ \\
  Position & $x$ & $\lambda$ \\
  Velocity & $v$ & $\bar v$ \\
  Distribution function& $f$ & $n_0/{\bar v}^3$ \\
  Electric field& $E$ & $\omega_p \bar v \frac{m}{q}$ \\
  Magnetic field& $B$ & $\omega_p \frac{m}{q}$ \\
  \hline
\end{tabular}
\caption{Units of Normalization. Here, $\bar v$ is the characteristic scale of $v$, $\omega_p=\sqrt{n_0 q^2/m\epsilon_0}$ is the electron plasma frequency with  total number of electrons by $n_0=\int_{\Omega_x\times\Omega_v} fd\mathbf{x}d\mathbf{v}$, and $\lambda=\bar v/\omega_p$.}
\label{Table:uni}
\end{table}
Choose the characteristic set of variables as they are shown in Table~\ref{Table:uni}. Set $\bar v=c$, we can derive the normalized equations,
\begin{equation}\label{eq:1}
\begin{aligned}
\frac{\partial f}{\partial t}+\mathbf{v}\cdot\frac{\partial f}{\partial\mathbf{x}}+(\mathbf{E}+\mathbf{v\times B})\cdot\frac{\partial f}{\partial\mathbf{v}}=0,\\
\frac{\partial \mathbf{E}}{\partial t}=\nabla\times \mathbf{B}-\int_{\Omega_v}\mathbf{v} fd\mathbf{v}, \quad \nabla\cdot \mathbf{B}=0,\\
\frac{\partial \mathbf{B}}{\partial t}=-\nabla\times \mathbf{E},\quad \nabla\cdot \mathbf{E}=\int_{\Omega_v} fd\mathbf{v}-\rho_0.
\end{aligned}
\end{equation}

In the magnetostatic limit, i.e. the case of $\frac{\partial \mathbf{B}}{\partial t}\rightarrow0$, the magnetic field can be taken as a background field. Denote the background magnetic field as $\mathbf{B}=\mathbf{B}_0(\mathbf{x})$, the charge density as $\rho(\mathbf{x},t)=\int_{\Omega_v}f d\mathbf{v}$ and the current density as $\mathbf{J}(\mathbf{x},t)=\int_{\Omega_v}\mathbf{v} fd\mathbf{v}$. The VM system (\ref{eq:1}) is  reduced to
\begin{align}
\frac{\partial f}{\partial t}+\mathbf{v}\cdot\frac{\partial f}{\partial\mathbf{x}}+ \left(\mathbf{E}+\mathbf{v}\times\mathbf{B}_0\right)\cdot\frac{\partial f}{\partial\mathbf{v}}=0,  \label{eq:2}\\
\frac{\partial\mathbf{E}}{\partial t}=\nabla\times\mathbf{B}_0-\mathbf{J},\quad  \nabla\cdot\mathbf{B}_0=0,  \label{eq:3}\\
\nabla\times\mathbf{E}=0,  \quad \nabla\cdot\mathbf{E}=\rho-\rho_0.   \label{eq:4}
\end{align}
The total charge neutrality of the reduced system is ensured due to $\int_{\Omega_x} (\rho-\rho_0)d\mathbf{x}=0$.

Notice that Maxwell's equations (3) and (4) can be decoupled, thus the system can be treated in two different ways.
The first way is to consider Amp{\`e}re's equation,
\begin{equation}\label{eq:AE}
\begin{gathered}
 \frac{\partial\mathbf{E}}{\partial t}=\nabla\times\mathbf{B}_0-\mathbf{J}.
\end{gathered}
\end{equation}
It together with Eq.(\ref{eq:2}) is called the Vlasov--Amp{\`e}re equation~\citep{Nicolas2011A,GChen2011,YingdaCheng2014,Perin2015Hamiltonian}.

Another way is to consider Faraday's law and Gauss's equation,
$$\nabla\times\mathbf{E}=0,\quad \nabla\cdot\mathbf{E}=\rho-\rho_0.$$
By introducing the electric potential $\phi_{f}$, it follows that
\begin{equation}
\begin{gathered}
{\bf E}(\mathbf{x},t)=-\nabla\phi_{f}(\mathbf{x},t),\quad -\Delta\phi_{f}= {\rho(\mathbf{x},t)-\rho_0}.\label{eq:poisson}
\end{gathered}
\end{equation}
Combining Eqs.(\ref{eq:2}) and (\ref{eq:poisson}) derives the so-called  Vlasov--Poisson equations.

As above, we have introduced two different reduced
equations:
the Vlasov--Poisson equation (\ref{eq:2},\ref{eq:poisson}) and the   Vlasov--Amp{\`e}re equation (\ref{eq:2},\ref{eq:AE}). Moreover, the two equations are connected.
By taking the divergence and curl of Amp{\`e}re's equation (\ref{eq:AE}), Eqs.(\ref{eq:4}) can be derived under the continuity condition $\frac{\partial \rho}{\partial t}+\nabla\cdot \mathbf{J}=0$ and $\nabla\times \mathbf{J}=\Delta \mathbf{B}_0$, and the following restriction of the initial electric field,
$$\nabla\cdot \mathbf{E}(\mathbf{x},0)-\rho(\mathbf{x},0)+\rho_0=0, \quad\nabla\times\mathbf{E}(\mathbf{x},0)=0.$$
Therefore, the solution of the Vlasov--Amp{\`e}re equation under the above  restrictions is a solution to the Vlasov--Poisson equation.
In one dimensional case, in order to connect the solution of the Vlasov--Poisson equation to the solution of the Vlasov--Amp{\`e}re equation~\citep{Back2014F}, we can add the requirement $\int_{\Omega_x} E(x,t)dx=0$ and the initial constraint $\int_{\Omega_x\times\Omega_v} vf(x,v,0) dxdv=0$.

In this paper, we mainly focus on the numerical study for the Vlasov--Poisson equations described by (\ref{eq:2}) and (\ref{eq:poisson}). In practical computation, the boundary conditions can be taken as zero, i.e., $f(\mathbf{x,v},t)=0$ and $\phi_f(\mathbf{x})=0$ on $\partial\Omega_x$. Also the periodic boundary condition is used frequently.

\section{Poisson bracket structure of Vlasov--Poisson equations}
In this section, we introduce the Poisson bracket whose fundamental definition can be found in ~\citep{Olver1986}. With the following defined Poisson bracket, in this section we describe  the Vlasov--Poisson equations (\ref{eq:2}) and (\ref{eq:poisson}) as the Poisson bracket system.

Suppose $\mathcal{F}$ and $\mathcal{G}$ are two functionals of $f$, we  define the bracket of two functionals as
\begin{equation}
\begin{aligned}
\{\{\mathcal{F},\mathcal{G}\}\} & (f)
=\int_{\Omega_{x}\times\Omega_{v}} f\left\{ \frac{\delta\mathcal{F}}{\delta f},\frac{\delta\mathcal{G}}{\delta f}\right\}_{\mathbf{xv}}d\mathbf{x}d\mathbf{v}
+\int_{\Omega_{x}\times\Omega_{v}} f{\bf B}_0\cdot\left(\frac{\partial}{\partial\mathbf{v}}\frac{\delta\mathcal{F}}{\delta f}\times\frac{\partial}{\partial\mathbf{v}}\frac{\delta\mathcal{G}}{\delta f}\right)d\mathbf{x}d\mathbf{v},
\end{aligned}
\label{eq:MMWB}
\end{equation}
where $\frac{\delta \mathcal{G}}{\delta f}$ is the variational derivative \footnote{The variational derivative~\citep{Morrison1998Hamiltonian,Marsden1999Introduction} $\frac{\delta \mathcal{G}}{\delta f}$ of a functional $\mathcal{G}[f]$ is defined by
\begin{align*}
\int_{\Omega_{x}\times\Omega_{v}}\frac{\delta\mathcal{G}}{\delta f}\delta f d\mathbf{x}d\mathbf{v}=\lim\limits_{\epsilon \to 0}\frac{\mathcal{G}[f+\epsilon\delta f]-\mathcal{G}[f]}{\epsilon}.
\end{align*}}, and the operator $\left\{ \cdot,\cdot\right\} _{\mathbf{xv}}$ is the canonical Poisson bracket of two functions $m(\mathbf{x,v})$ and  $n(\mathbf{x,v})$, that is
\[
\{m,n\}_{\mathbf{xv}}=\frac{\partial m}{\partial \mathbf{x}}\cdot\frac{\partial n}{\partial \mathbf{v}}-\frac{\partial m}{\partial \mathbf{v}}\cdot\frac{\partial n}{\partial \mathbf{x}}.
\]
It can be checked easily that the bracket (\ref{eq:MMWB}) is bilinear and anti-symmetric. Furthermore, when $\mathbf{B}_0$ satisfies $\nabla\cdot\mathbf{B}_{0}(\mathbf{x})=0$, the bracket yields the Jacobi identity which has been proved in \ref{app:1}.

With the bracket (\ref{eq:MMWB}), we can present the following system,
\begin{equation}
\frac{d\mathcal{F}}{d t}=\{\{\mathcal{F},\mathcal{H}\}\},\label{eq:PoissonVM}
\end{equation}
where $\mathcal{F}$ is any functional of the solution $f$, and $\mathcal{H}$ is the Hamiltonian functional and the global energy of the system,
\begin{equation}
\begin{aligned}
\mathcal{H}[f]
& =\frac{1}{2}\int_{\Omega_{x}\times\Omega_{v}}\mathbf{v}{}^{2}fd\mathbf{x}d\mathbf{v}+\frac{1}{2}\int_{\Omega x}\mathbf{E}^{2}d\mathbf{x}\\
& =\frac{1}{2}\int_{\Omega_{x}\times\Omega_{v}}\mathbf{v}{}^{2}fd\mathbf{x}d\mathbf{v}+\frac{1}{2}\int_{\Omega_{x}\times\Omega_{v}}\phi_{f}fd\mathbf{x}d{\bf v}-\frac{1}{2}\rho_0\int_{\Omega_{x}}\phi_{f}d\mathbf{x}.\label{eq:HamiltonVM}
\end{aligned}
\end{equation}
In fact, by setting $\mathcal{F}[f]=\int_{\Omega_x\times\Omega_v}f(\tilde{\mathbf{x}},\tilde{\mathbf{v}},t)\delta(\mathbf{x}-\tilde{\bf x})\delta(\mathbf{v}-\tilde{\bf v})d\tilde{\bf x}d\tilde{\bf v}$, and defining the local energy $h(\mathbf{x},\mathbf{v})=\frac{\delta\mathcal{H}}{\delta f}(\mathbf{x},\mathbf{v})={\bf v^{2}}/2+\phi_{f}({\bf x})$, the formulation (\ref{eq:PoissonVM}) leads to
\begin{align}\label{eq:poisson3}
\frac{\partial f}{\partial t}&=-\{f,h\}_{\mathbf{xv}}-\mathbf{B}_0\cdot\left(\frac{\partial f}{\partial \mathbf{v}}\times \frac{\partial h}{\partial \bf v}\right).
\end{align}
Eq. (\ref{eq:poisson3}) recovers the Vlasov equation (\ref{eq:2}).

As a reduced model from the Vlasov--Maxwell equations, conservative properties of the Vlasov--Poisson system is also of great importance. With the help of the Poisson bracket (\ref{eq:MMWB}) and the formulation (\ref{eq:PoissonVM}), the conservative properties can be redescribed as follows.

\noindent{\bf Energy conservation.}
The total energy defined in (\ref{eq:HamiltonVM}) is invariant as
\[\frac{d\mathcal {H}}{dt}=\{\{\mathcal{H},\mathcal{H}\}\}=0.
\]

\noindent{\bf Momentum conservation.}
The total momentum $\mathcal{P}=\int_{\Omega_x\times\Omega_v} \mathbf{v}fd\mathbf{x} d\mathbf{v}$ is invariant when $\mathbf{B}_0=\mathbf{0}$.

In fact, taking the derivative of $\mathcal{P}$
gives
\begin{equation*}
\begin{aligned}
\frac{d\mathcal {P}}{dt}&=\{\{\mathcal{P},\mathcal{H}\}\}
&=-\int_{\Omega_x} \nabla\phi_f \rho d\mathbf{x}+\int_{\Omega_x\times\Omega_v} (\mathbf{v}\times \mathbf{B}_0)f d\mathbf{x}d\mathbf{v}.
\end{aligned}
\end{equation*}
The first term vanishes  due to the Poisson equation which is expressed as $\nabla\cdot \mathbf{E}=\rho-\rho_0, \nabla\times \mathbf{E}=0$, because
\begin{align*}
-\int_{\Omega_x} \nabla\phi_f \rho d\mathbf{x}
&=\int_{\Omega_x} \mathbf{E}(\nabla\cdot \mathbf{E}+\rho_0) d\mathbf{x}\\
&=-\int_{\Omega_x} (\nabla \times \mathbf{E})\times \mathbf{E} d\mathbf{x}+\rho_0\int_{\Omega_x} \mathbf{E}d\mathbf{x}=0.
\end{align*}
Here, $\int_{\Omega_x} \mathbf{E}d\mathbf{x}=0$ for zero or periodic boundary condition of $\phi_f$.

\noindent{\bf Casimir functional.}
Any functional in the form $\mathcal{C}[f]=\int_{\Omega_x\times\Omega_v}C(f)d\mathbf{x}d\mathbf{v}$ is invariant, as it is
a Casimir functional~\citep{Morrison1998Hamiltonian} w.r.t the Poisson bracket (\ref{eq:MMWB}) with $\mathbf{B}_0=\mathbf{0}$.
In fact, taking the bracket operation between $\mathcal{C}[f]$ and any functional $\mathcal{G}[f]$ gives
\[\{\{\mathcal{C},\mathcal{G}\}\}=\int_{\Omega_x\times\Omega_v} \left\{f,\frac{\delta\mathcal{C}}{\delta f}\right\}_{\mathbf{xv}}\frac{\delta\mathcal{G}}{\delta f} d{\mathbf{x}}d{\mathbf{v}}.
\]
The above equation vanishes because $\frac{\delta \mathcal{C}[f]}{\delta f}=C'(f)$ is a function of $f$. Therefore, $\mathcal{C}[f]$ is an invariant for any PDE of $f$  equipped with the Poisson bracket (\ref{eq:MMWB}), including the Vlasov equation, that is
\begin{equation*}
\frac{d\mathcal{C}}{dt}=\{\{\mathcal{C},\mathcal{H}\}\}=0.
\end{equation*}

This provides a class of conservative quantities of the Vlasov--Poisson system~(\ref{eq:2},\ref{eq:poisson}). For instance,  the integral norm $I=\int_{\Omega_x\times\Omega_v} f^p d\mathbf{x} d\mathbf{v}$ and the entropy $S=\int_{\Omega_x\times\Omega_v} f \ln f d\mathbf{x}d\mathbf{v}$ of the system are invariant.

\section{Spatial discretization for Vlasov--Poisson equations}
In this section, we present the construction of numerical approach for solving the Vlasov--Poisson equations (\ref{eq:2}) and (\ref{eq:poisson}).

According to the idea of PIC, the distribution function is sampled by a set of super particles. Assume that the distribution function $f$ is approximated by the Klimontovich representation
\begin{equation}
f_{h}(\mathbf{x},\mathbf{v},t)=\sum_{s=1}^{N_p}\omega_{s}\delta(\mathbf{x}-\mathbf{X}_{s}(t))\delta(\mathbf{v}-\mathbf{V}_{s}(t)),\label{eq:Disf}
\end{equation}
where $(\mathbf{X}_{s},\mathbf{V}_{s})$ is the s-th super particle's coordinates in phase space, and $\omega_{s}$ is the particle weight, $N_p$ is the number of super particles.
The Vlasov equation (\ref{eq:2}) is then transformed to particle equations,
\begin{gather}
\dot{\mathbf{X}}_{s}=\mathbf{V}_{s},\quad
\dot{\mathbf{V}}_{s}=\int_{\Omega_x}\left({\bf E}(\mathbf{x},t)+{\bf V}_s\times {\bf B}_0(\mathbf{x})\right)\delta(\mathbf{x}-{\bf X}_s)d\mathbf{x}, \quad
 s=1,2,\ldots N_p.
\label{eq:Vlasov-1}
\end{gather}

We use finite element methods to solve the Poisson equation (\ref{eq:poisson}) with the Dirichlet boundary condition, i.e.
\begin{gather*}
-\Delta\phi_f =\rho-\rho_0 \text{~in~}\Omega_{x},\qquad
\phi_f(x,t)=0\text{~on~}\partial\Omega_{x}.
\end{gather*}
Let $V$ be a Hilbert space and $V'$ be the dual space of it.
In the framework of finite element, the variational problem of Poisson equation (\ref{eq:poisson}) is to find solutions $\phi\in V$ such that for any $\psi\in V$ the following equation holds
\begin{equation}
\left(\nabla\phi,\mathbf{\nabla}\psi\right)=\left<\rho-\rho_0,\psi\right>,\quad\forall \psi\in V.
\label{eq:VariForm}
\end{equation}
Here, $\rho-\rho_0\in V'$, $(\,\cdot, \,\cdot)$ is the inner product of $V$ and $\left<\,\cdot, \,\cdot\right>$ is the pairing of elements of $V'$ and $V$. It is known by the Lax-Milgram theorem that the solution to the variational problem (\ref{eq:VariForm}) is unique~\citep{Lions1967,Brenner2008}.

In order to determine the proper space $V$, we refer to ~\citep{HY2016,Kraus2017,Martin2021Variational} where finite element discretizations have been presented  to derive Hamiltonian algorithms for the Vlasov--Maxwell equations. Given a linear operator $D\in\mathrm\{grad,curl,div\}$, we denote the Sobolev spaces as
\begin{align*}
&H(D,\Omega):=\{v\in L^2(\Omega),Dv\in L^2(\Omega)\},\\
&H_0(D,\Omega):=\{v\in H(D,\Omega), \text{trace~} v=0 \text{~on~} \partial \Omega\}.
\end{align*}
Denote $H^h(D,\Omega)$ as the finite element subspace of $H(D,\Omega)$.
For the finite element PIC discretisation of the Vlasov--Maxwell equations (\ref{eq:1}), we have the following results. Here the perfect conducting boundary (PEC) conditions are considered.
\begin{pro}[~\citep{HY2016}]
Suppose that the fields ${\bf E}$ and ${\bf B}$ are discretised in spaces ${\mathcal{E}}_h\subset H_0(curl,\Omega)$ and ${\mathcal{B}}_h\subset H_0(div,\Omega)$ respectively. The approximate problem for Maxwell's equations is to find $(\mathbf{E}_h, \mathbf{B}_h)\in \mathcal{E}_h\times\mathcal{B}_h$ such that
\begin{equation*}
\begin{gathered}
\left(\partial_{t}\mathbf{E}_h,\mathbf{\Phi}\right)= \left(\mathbf{B}_h,\nabla\times\mathbf{\Phi}\right)- \left(\sum_{s}\omega_s\mathbf{V}_{s}\delta(\mathbf{x}-\mathbf{X}_{s}),\mathbf{\Phi}\right), \quad\forall \mathbf{\Phi}\in\mathcal{E}_h,\\
\left(\partial_{t}\mathbf{B}_h,\mathbf{\Psi}\right)=-\left(\nabla\times\mathbf{E}_h,\mathbf{\Psi}\right),\quad \forall \mathbf{\Psi}\in \mathcal{B}_h.
\end{gathered}
\end{equation*}
If $\nabla\times {\mathcal{E}}_h\subset {\mathcal{B}}_h$ and $\nabla\cdot{\bf{B}}_h=0$, the resulting
semi-discrete system is Hamiltonian with a discrete Poisson structure.
\end{pro}

To understand the above theorem, we can  employ the following diagram. That is,  in order to derive a discrete Hamiltonian system, we need choose the finite element to satisfies the de Rham Complex diagram
\begin{center}
 $ \CD
   H(grad) @>grad >> H(curl) @>curl>> H(div) @>div>> L^2 \\
   @V \pi VV @V \pi VV @V \pi VV @V \pi VV  \\
   H^h(grad) @>grad>> H^h(curl) @>curl>>H^h(div) @>div>> L^{2,h}
 \endCD $
\end{center}

According to the above comment and the relation $\mathbf{E}=-\nabla \phi$,
for $\mathbf{E}\in H_0(curl,\Omega)$  we can choose $\phi\in V=H_0(grad,\Omega_x)$.

The approximate problem of (\ref{eq:VariForm}) is then to find $\phi_{h}\in V_{h}\subset H_{0}^{1}(\Omega_{x})$
such that
\begin{equation}
\left(\nabla\phi_{h},\mathbf{\nabla}\psi_{h}\right)=\left<\rho_h,\psi_{h}\right>,\quad\forall \psi_{h}\in V_{h}.
\label{eq:DisVariForm}
\end{equation}
Here, $\rho_h\in L^{2}(\Omega_{x})$ is the discrete version of $\rho-\rho_0$. The inner product is defined by $\left(\mathbf{f},\mathbf{g}\right)=\int_{\Omega_{x}}\mathbf{f}\cdot\mathbf{g}d\mathbf{x}$ on space $L^{2}(\Omega_{x})$ and $V_{h}$ is the finite-dimensional subspace. In particular, for all $\psi_{h}\in H_{0}^{1}(\Omega_{x})$ and $\rho_h\in L^{2}(\Omega_{x})\subset H^{-1}(\Omega_{x})$, we have $\left(\rho_h,\psi_{h}\right)=\left<\rho_h,\psi_{h}\right>$,
and the solution to the variational problem (\ref{eq:DisVariForm}) is unique.
Thus, the well studied finite element spaces including linear element, quadratic Lagrange element, and Spline element etc. can be implemented.

Suppose that $\{W_j(\mathbf{x})\}_{j=1}^{N}$ are piecewise polynomial basis functions of the space $V_{h}$, then $\mathbf{\phi}_{h}\in V_{h}$ can be expressed as
\begin{equation}
\mathbf{\phi}_{h}(\mathbf{x},t)=\sum_{j=1}^{N}\phi_{j}(t){W}_{j}(\mathbf{x}).\label{eq:DisEB}
\end{equation}

In order to get an appropriate regularization of $\rho_h$, we introduce the regularizing function $S_{\epsilon}(\mathbf{x})$ which satisfies $\int_{\mathbb{R}^d}S(\mathbf{x})d\mathbf{x}=1$ and $S_{\epsilon}(\mathbf{x})=\frac{1}{\epsilon^d}S(\frac{\mathbf{x}}{\epsilon})$~\citep{Raviart1985An,Cottet1984Particle,Victory1989On}. Here $\epsilon$ denotes the width of the kernel and $d=1,2,3$.

Then, we define $f^{\epsilon}_h=f_h\ast S_{\epsilon}$, that is $f^{\epsilon}_{h}(\mathbf{x},\mathbf{v},t)=\sum\limits_{s=1}^{N_p}\omega_{s}S_{\epsilon}(\mathbf{x}-\mathbf{X}_{s}(t))\delta(\mathbf{v}-\mathbf{V}_{s}(t))$.
By means of $S_{\epsilon}$,  we can choose $\rho_h(\mathbf{x})=\sum\limits_{s=1}^{N_p}\omega_{s}S_{\epsilon}(\mathbf{x}-\mathbf{X}_{s})-\rho_0$ derived from $\rho_h(\mathbf{x},t)=\int_{\Omega_v}f^{\epsilon}_h(\mathbf{x},\mathbf{v},t)d\mathbf{v}-\rho_0$. On the other hand, $\rho_h$ can be projected to the space $V_h$, which gives $\rho_h(\mathbf{x})=\sum\limits_{i=1}^N\rho_{j}{W}_{j}(\mathbf{x})$. This   will help us in computation due to that $({W}_{i},{W}_{j})$ can be calculated easily, and $\rho_h\in V_h\subset H_0^1\subset L^2$.

Substitute Eq.(\ref{eq:DisEB}) into Eq.(\ref{eq:DisVariForm})
and take $\psi_h={W}_{i}$, then we get the following approximate
problem in  matrix formulation,
\begin{equation}
\sum_{j=1}^{N}\left(\nabla{W}_{j},\nabla{W}_{i}\right)\phi_{j}=\left(\rho_h,{W}_{i}\right),i=1\ldots N,
\label{eq:semiDis}
\end{equation}
where the solutions $\{\phi_j\}_{i=1}^N$ are the  functions of particle positions $\mathbf{X}_s$, $s=1,2,\ldots,N_p$. In order to get a vector expression for the particles, we denote $\mathbf{X}=({\bf X}_1^{\rm T},{\bf X}_2^{\rm T},\ldots,{\bf X}_{N_p}^{\rm T})^\mathrm{T} \in \mathbb{R}^{3N_p}$ and $\mathbf{V}=({\bf V}_1^{\rm T},{\bf V}_2^{\rm T},\ldots,{\bf V}_{N_p}^{\rm T})^\mathrm{T}\in \mathbb{R}^{3N_p}$.
Then the approximate electric field is
\begin{equation}\label{eq:Eh}
\begin{aligned}
&{\bf E}_{h}({\bf x},{\bf X},t)=-\sum_{j=1}^{N}\phi_{j}({\bf X}(t))\nabla {W}_{j}(\mathbf{x}).
\end{aligned}
\end{equation}
Substituting $\mathbf{E}_h$ into Eq. (\ref{eq:Vlasov-1}), we get the discrete particle equations for $s=1,2,\ldots N_p$,
\begin{equation}
\begin{aligned}
\dot{\mathbf{X}}_{s}&=\mathbf{V}_{s},\\
\dot{\mathbf{V}}_{s}&=-\sum_{j=1}^{N}\phi_{j}({\bf X})\nabla {W}_{j}(\mathbf{X}_s)+{\bf V}_s\times {\bf B}_0(\mathbf{X}_s).
\end{aligned} \label{eq:semiDisV}
\end{equation}

Define the discrete bracket by discretizing the Poisson bracket (\ref{eq:MMWB}) shown in \ref{app:2}. Moreover, it is proved in \ref{app:3} that the following discrete bracket is Poisson as long as $\nabla\cdot {\bf B}_0({\bf x})=0$. With the defined discrete Poisson bracket, the semi-discrete system (\ref{eq:semiDisV}) is a Hamiltonian ODE system. For two arbitrary functions $F$ and $G$ of $(\mathbf{X,V})$, it reads
\begin{equation}
\begin{aligned}
\left\{ F,G\right\} (\mathbf{X},\mathbf{V})
=\sum_{s=1}^{N_p}\frac{1}{\omega_{s}}\left(\frac{\partial F}{\partial\mathbf{X}_{s}}\cdot\frac{\partial G}{\partial\mathbf{V}_{s}}-\frac{\partial G}{\partial\mathbf{X}_{s}}\cdot\frac{\partial F}{\partial\mathbf{V}_{s}}\right)+\sum_{s=1}^{N_p} \frac{1}{\omega_{s}}{\bf B}_0(\mathbf{X}_s)\cdot\left(\frac{\partial F}{\partial\mathbf{V}_s} \times\frac{\partial G}{\partial\mathbf{V}_s}\right).
\end{aligned}
\label{eq:dPoisson}
\end{equation}

Accordingly, by inserting (\ref{eq:Disf}) into (\ref{eq:HamiltonVM}) we can derive the following discrete Hamiltonian function,
\begin{equation}
\begin{aligned}
H({\bf   X,V})=\mathcal{H}[f_h]
& =\frac{1}{2}\int_{\Omega_{x}\times\Omega_{v}}\mathbf{v}{}^{2}f_h d\mathbf{x}d\mathbf{v}+\frac{1}{2}\int_{\Omega_x}(\nabla\phi_h)\cdot (\nabla\phi_h) d\mathbf{x}\\
& =\frac{1}{2} \sum_{s=1}^{N_p} \omega_s\mathbf{V}_s^{2}+\frac{1}{2}\sum_{j=1}^N\sum_{k=1}^N \phi_j({\bf X})\phi_k({\bf X})\int_{\Omega_x}\nabla W_j(\mathbf{x})\cdot\nabla W_k(\mathbf{x}) d\mathbf{x}
\end{aligned}
\label{eq:dH}
\end{equation}
For the finite element discretization, the discrete Hamiltonian function~(\ref{eq:dH}) is natural, and can be  rewritten in matrix form easily. With the discrete Poisson bracket~(\ref{eq:dPoisson}) and the discrete Hamiltonian~(\ref{eq:dH}), the semi-discrete equations of motion (\ref{eq:semiDisV}) can be recovered by
\begin{equation}
\dot{\bf X}_s=\{{\bf X}_s,H\}, \quad \dot{\bf V}_s=\{{\bf V}_s,H\}.
\label{eq:sdeq}
\end{equation}

As follows, we express the semi-discretised equations~(\ref{eq:sdeq}) by their matrix formulation.
Introduce the following notations for the field basis and variables:
\begin{equation*}
\begin{gathered}
\mathbf{\Phi}({\bf X})=(\phi_1,\phi_2,\ldots,\phi_{N})^{\rm T}({\bf X}) \in \mathbb{R}^{N},\quad \mathbb{M}\in \mathbb{R}^{N\times N},\text{~with~} \mathbb{M}_{jk}=\int_{\Omega x}\nabla W_j({\bf x})\cdot\nabla W_k({\bf x}) d\mathbf{x},\\
\mathbb{G}({\bf X})=\left(
          \begin{array}{cccc}
            \nabla W_1({\bf X}_1) & \nabla W_2({\bf X}_1) & \cdots & \nabla W_N({\bf X}_1) \\
            \nabla W_1({\bf X}_2) & \nabla W_2({\bf X}_2) & \cdots & \nabla W_N({\bf X}_2) \\
            \cdots & \cdots & \cdots &\cdots \\
            \nabla W_1({\bf X}_{N_p}) & \nabla W_2({\bf X}_{N_p}) & \cdots & \nabla W_N({\bf X}_{N_p}) \\
          \end{array}
        \right)\in \mathbb{R}^{(3N_P)\times N},\\
\hat{\mathbb{B}}(\mathbf{X})=\rm{diag}(\hat{\bf B}_0({\bf X}_1),\hat{\bf B}_0({\bf X}_2),\ldots,\hat{\bf B}_0({\bf X}_{N_p}))\in \mathbb{R}^{(3N_p)\times(3N_p)}.
\end{gathered}
\end{equation*}
where $\hat{\bf B}_0$ is the hat matrix w.r.t the vector ${\bf B}_0$. The hat matrix of a vector ${\bf B}=(B_x,B_y,B_z)$ is defined as
\begin{equation}
\hat{\bf B}=\left(
            \begin{array}{ccc}
              0 & B_z & -B_y \\
              -B_z & 0 & B_x \\
              B_y & -B_x & 0 \\
            \end{array}
          \right). \label{eq:hatB}
\end{equation}
Denote the diagonal weight matrix $\mathbf{\Omega}=\mathrm{diag}(\omega_1,\omega_2,\ldots,\omega_{N_p})\in \mathbb{R}^{N_p\times N_p}$, and
$\mathbf{I}$ as the $3$-dimensional identity matrix. Let $$\mathbb{N}=\mathbf{\Omega}\otimes \mathbf{I}\in \mathbb{R}^{(3N_p)\times (3N_p)}.$$
Therefore, the discrete Poisson bracket~(\ref{eq:dPoisson}) can be reformulated in matrix form as
\begin{equation}
\begin{aligned}
\left\{ F,G\right\}(\mathbf{Z}) &=\sum_{s=1}^{N_p} \frac{1}{\omega_{s}}\left(
                                                \begin{array}{cc}
                                                \left(\frac{\partial F}{\partial\mathbf{X}_s}\right)^{\mathrm{T}}, & \left(\frac{\partial F}{\partial\mathbf{V}_s}\right)^{\rm T} \\
                                                \end{array}
                                              \right)
                                              \left(
                                                \begin{array}{cc}
                                                  0 & \mathbf{I} \\
                                                  -\mathbf{I} & \hat{\mathbf{B}}_0(\mathbf{X}_s) \\
                                                \end{array}
                                              \right)
                                              \left(
                                                \begin{array}{c}
                                                  \frac{\partial G}{\partial\mathbf{X}_s} \\
                                                  \frac{\partial G}{\partial\mathbf{V}_s} \\
                                                \end{array}
                                              \right)\\
&=\frac{\partial F}{\partial\mathbf{Z}}^T \mathbb{K}({\bf X}) \frac{\partial G}{\partial\mathbf{Z}},
\end{aligned}
\end{equation}
where $\mathbf{Z}=(\mathbf{X},\mathbf{V})$, and $\mathbb{K}({\bf X})=\left(
                \begin{array}{cc}
                  0 & \mathbb{N}^{-1} \\
                  -\mathbb{N}^{-1} & \mathbb{N}^{-1}\hat{\mathbb{B}}({\bf X}) \\
                \end{array}
              \right)$ is the Poisson matrix.
With the above notations, the discrete Hamiltonian (\ref{eq:dH}) can be rewritten as,
\begin{equation}
\begin{aligned}
H(\mathbf{X,V})=\frac{1}{2} \mathbf{V}^T \mathbb{N}\mathbf{V}+\frac{1}{2}\mathbf{\Phi(X)}^{\mathrm T} \mathbb{M} \mathbf{\Phi(X)}.
\end{aligned}
\end{equation}
The corresponding system (\ref{eq:semiDisV}) can be rewritten with Poisson matrix $\mathbb{K}(\mathbf{X})$, which is
\begin{equation}
\begin{aligned}
\dot {\bf X}&={\bf V},\\
\dot{\bf V}&=\mathbb{G}({\bf X})\mathbf{\Phi(X)}+ \hat{\mathbb{B}}({\bf X}){\bf V}.
\end{aligned}
\end{equation}

With the help of the discrete Poisson bracket, the conservation of the discrete energy $H({\bf X,V})$ in Eq.(\ref{eq:dH}) can be easily verified, that is
\[
\frac{dH}{dt}=\{H,H\}(\mathbf{X},\mathbf{V})=0.
\]
Similarly, the total charge $C=\int_{\Omega_x\times\Omega_v}f_h d\mathbf{x}d\mathbf{v}=\sum_{s=1}^{N_p}\omega_s$ is invariant, as it is a Casimir function of the discrete Poisson bracket (\ref{eq:dPoisson}),
$$\frac{dC}{dt}=\{C,H\}(\mathbf{X},\mathbf{V})\equiv 0.$$

%

\section{Temporal discretization for Vlasov--Poisson equations}
In this section we are designing temporal discretisation for the Vlasov--Poisson equations that can preserve the Poisson bracket structure. As the Poisson bracket is not canonical, traditional time integrators such as Runge-Kutta methods can not be used to construct Poisson-structure-preserving methods. The idea of splitting technique is to split the original system into several subsystems. Each subsystem can be solved exactly, and has the same structure as ones for the original system~\citep{Yoshida1990,McLachlan2002}. Thus, the composition of solutions to subsystems leads to numerical methods which can preserve the structure of the original systems. High-order Hamiltonian splitting methods have been developed for Vlasov--Maxwell equations~\citep{He2015Hamiltonian,LYZ2019} and Vlasov--Poisson equations without magnetic fields~\citep{Nicolas2015H,Casas2017}.
In this section, for the non-canonical Hamiltonian ODE system (\ref{eq:semiDisV}) we establish the Poisson-structure-preserving methods via Hamiltonian splitting.

We split the Hamiltonian as,
\begin{align*}
 & H=H_v+H_e, \quad
  \text{where~} H_v=\frac{1}{2} \sum_{s=1}^{N_p} \omega_s\mathbf{V}_s^{2}, \quad H_e=\frac{1}{2}\mathbf{\Phi}(\mathbf{X})\mathbb{M}\mathbf{\Phi}(\mathbf{X}).
\end{align*}
With each part of the Hamiltonian, we can split the system into two parts,
$$\dot{\mathbf{Z}}=\{\mathbf{Z},H_v\},\quad \dot{\mathbf{Z}}=\{\mathbf{Z},H_e\}.$$
Each subsystem possesses the same Poisson bracket structure as the system (\ref{eq:semiDisV}). Thus the composition of solutions to each subsystem preserves the Poisson bracket of the original system due the group property of Poisson bracket structure.

Associated with the Hamiltonian
$H_{v}$, the subsystem is $\dot{F}=\left\{ F,H_{v}\right\} $. For the s-th particle it is,
\begin{equation}
\begin{aligned}
& \mathbf{\dot{X}}_{s}=\mathbf{V}_{s},\\
 & \dot{\mathbf{V}}_{s}=\mathbf{V}_{s}\times \mathbf{B}_0(\mathbf{X}_{s}).
\end{aligned}
\label{eq:dSysHE}
\end{equation}

When $\mathbf{B}_0$ is constant and $|\mathbf{B}_0|=b$, the exact update mapping of this subsystem~(\ref{eq:dSysHE}) with step
size $\Delta t$ is
\begin{equation}
\phi^{H_v}(\Delta t):
\begin{aligned}
& \mathbf{X}_{s}(t+\Delta{t})=\mathbf{X}_{s}(t)+(\Delta{t}\mathbf{I}+\frac{1-\cos{b\Delta{t}}}{b^2}\hat{\mathbf{B}}_0+\frac{b\Delta{t}-\sin{b\Delta{t}}}{b^3}{\hat{\mathbf{B}}_0}^2)\mathbf{V}_{s}(t),\\
 & \mathbf{V}_{s}(t+\Delta{t})=(\mathbf{I}+\frac{\sin{b\Delta{t}}}{b}\hat{\mathbf{B}}_0+\frac{1-\cos{b\Delta{t}}}{b^2}{\hat{\mathbf{B}}_0}^2)\mathbf{V}_{s}(t).
\end{aligned}
\label{eq:dSolHE}
\end{equation}
where $\hat{\mathbf{B}}_0$ is the hat matrix as in Eq.(\ref{eq:hatB}). Following the idea, we  can also consider the problem in which  $\mathbf{B}_0$ is  non-homogeneous. In this case, we need to split the subsystem corresponding to $H_v$  further according to its new decomposition $H_v=\frac{1}{2}\omega_s\mathbf{V}_{sx}^{2}+\frac{1}{2}\omega_s\mathbf{V}_{sy}^{2}+\frac{1}{2}\omega_s\mathbf{V}_{sz}^{2}$.
More detail can be seen in~\citep{He2015Hamiltonian}.

The equation $\dot{F}=\left\{ F,H_{e}\right\} $ associated with the Hamiltonian $H_{e}$ is
\begin{equation}
\begin{aligned} & \mathbf{\dot{X}}_{s}=0,\\
 & \dot{\mathbf{V}}_{s}=-\sum_{j=1}^{N}\phi_{j}(\mathbf{X}(t))\int_{\Omega_x}\nabla{W}_{j}(\mathbf{x})\delta(\mathbf{x}-\mathbf{X}_{s})d\mathbf{x}.
\end{aligned}
\label{eq:dSysHB}
\end{equation}
The exact update of this subsystem~(\ref{eq:dSysHB})  with step size
$\Delta t$ is
\begin{equation}
\phi^{He}(\Delta t):\begin{aligned} & \mathbf{X}_{s}(t+\Delta t)=\mathbf{X}_{s}(t),\\
 & \mathbf{V}_{s}(t+\Delta t)=\mathbf{V}_{s}(t)-\Delta t \sum_{j=1}^{N}\phi_{j}(\mathbf{X}(t))\int_{\Omega_x}\nabla{W}_{j}(\mathbf{x})\delta(\mathbf{x}-\mathbf{X}_{s})d\mathbf{x}.
\end{aligned}
\label{eq:dSolHB}
\end{equation}

With the given exact solutions to the subsystems,
Poisson integrators can be derived by the compositions of the sub-flows~(\ref{eq:dSolHE}) and (\ref{eq:dSolHB}).
For example, the Poisson method  of first order can be constructed by
\[
\Phi(\Delta t)=\phi^{He}(\Delta t)\circ\phi^{Hv}(\Delta t),
\]
and a second order symmetric method can be derived from
\begin{align*}
\Phi(\Delta t)= & \phi^{Hv}(\Delta t/2)\circ\phi^{He}(\Delta t)\circ\phi^{Hv}(\Delta t/2).
\end{align*}
From the above, it is clear that these methods constructed by splitting  can  be implemented easily. More important, these methods can preserve the discrete Poisson bracket.

At the end of this section, we present the algorithm framework which helps to understand the computation procedure.
\begin{algorithm}[ht]
\caption{Algorithm framework}
\hspace*{0.02in} {\bf Input:}
$f(\mathbf{x},\mathbf{v},t=0)$
\begin{algorithmic}[1]
\State Approximate the initial condition $f_0$ by the distribution $f^{\epsilon}_h(\mathbf{x},\mathbf{v},t=0)=\sum_{s=1}^{N_p}\omega_{s}S_{\epsilon}(\mathbf{x}-\mathbf{X}_{s}(0))\delta(\mathbf{v}-\mathbf{V}_{s}(0))$ where $(\mathbf{X}_{s}(0),\mathbf{V}_{s}(0))$ is a set of particles distributed in the phase space
according to the density function $f_0$.
\State Generate finite element mesh and assemble stiffness matrix $\mathbb{M}$.
\State Compute the density $\rho_h=\sum_{s=1}^{N_p}\omega_{s}S_{\epsilon}(\mathbf{x}-\mathbf{X}_{s})-\rho_0$.
\State Assemble  the load matrix $F$ according to the density $\rho_h$. That is $F_i=(\rho_h,{W}_{i})$, $i=1\ldots N$.
\State Solve the large sparse linear system $\mathbb{M}{\phi_h}=F$. (A lot of linear solvers can be used such as Conjugate Gradient.)
\State Interpolate the electric field at the particle position by using discrete electric potential ${\phi_h}$. We can just use the information of $\nabla{W}_{i}(\mathbf{x})$, $i=1\ldots N$.
\State Update the position and velocity of particles by solving particle equation. High order symplectic algorithms can be considered.
\State Repeat steps 3-7 till the final time T.
\end{algorithmic}
\hspace*{0.02in} {\bf Output:}
$f^{\epsilon}_h(\mathbf{x},\mathbf{v},t=T)$ or other required data\\
\end{algorithm}

\section{Numerical experiments}
 In this section, we present  numerical experiments by using numerical methods presented in the previous section. The resulting numerical discretizations are derived by combining splitting method in time and finite element discretization in space. To derive the numerical results more efficiently we employ the parallel technique which is run on high performance computing workstation of LSEC Lab.

\subsection{Test problems in 1+1-dimensional phase space}
In this case, we perform various problems: Landau damping, Two-stream stability and Bump-on-tail instability. In the numerical simulation for all problems, we use periodic boundary condition. The number of particles is chosen as $N_p=10^6$. The system considered here is modelled by Vlasov equation coupled with normalized Poisson equation
\begin{equation*}
\mathbf{E}=-\partial{\phi},\quad -{\partial}^2{\phi}=\rho-1.
\end{equation*}

{\bf Landau damping.} Landau damping is referred to as the damping of a collective mode of oscillations in plasmas without
collisions of charged particles.
It, commonly believed, is caused by the energy exchange between electromagnetic wave and particles. Landau damping is also a very popular benchmark problem for testing the numerical methods applied to the Vlasov--Poisson equation due to that there are many analytical results
in literature~\citep{Filbet2003,Nicolas2009}. Here, we consider the case of one spatial dimension. In this case,  we take the initial distribution function as
\begin{equation*}
f_0(x,v)=\frac{1}{\sqrt{2\pi}}\exp(-\frac{v^2}{2})(1+\alpha\cos(kx)),
\end{equation*}
where the perturbation parameter $\alpha=0.001$ and $k$ is the wave number.

\begin{figure}[h!]
\centering
\subfigure[]
{
\includegraphics[scale=.5]{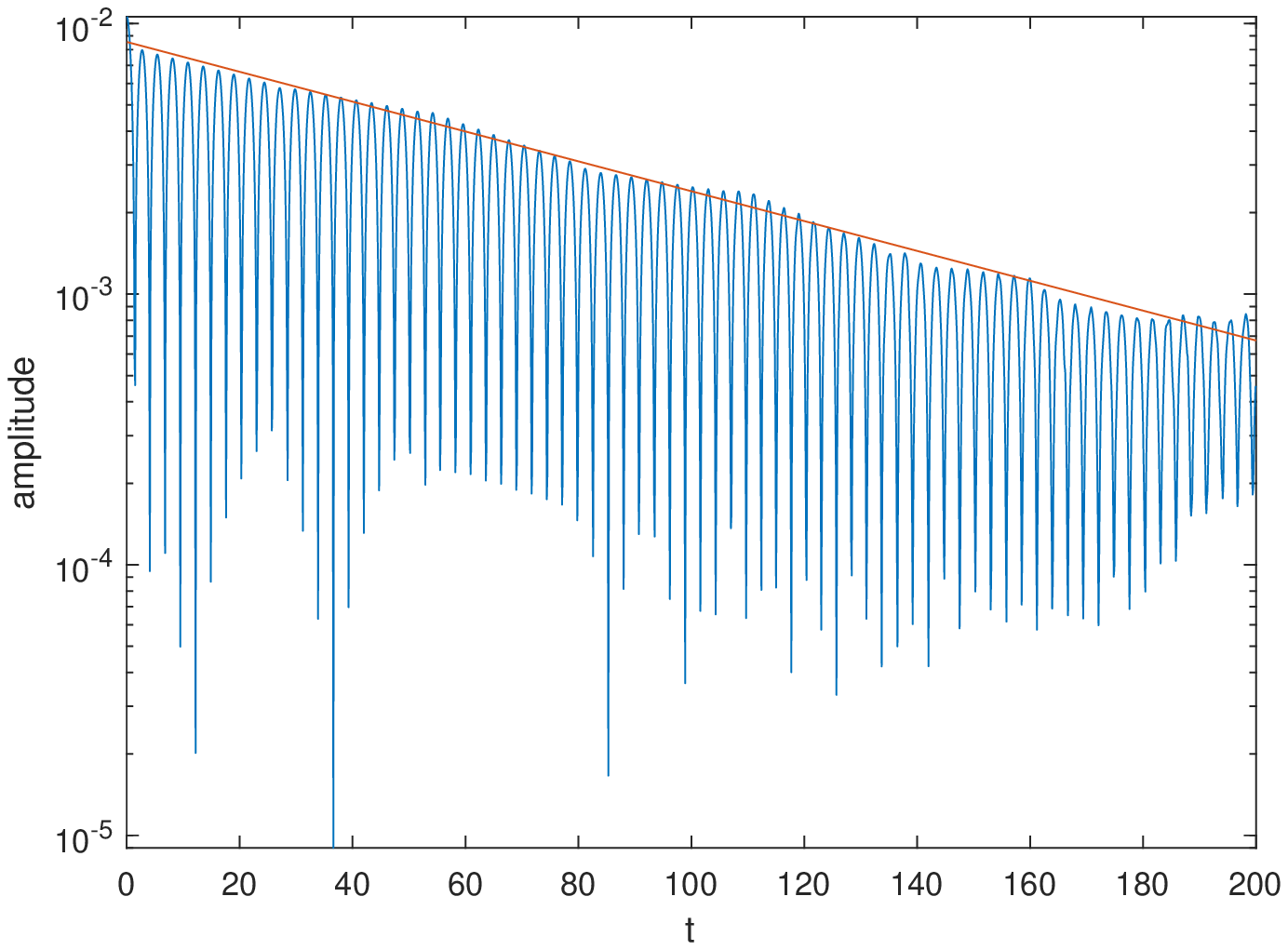}
}
\quad
\subfigure[]{
\includegraphics[scale=.5]{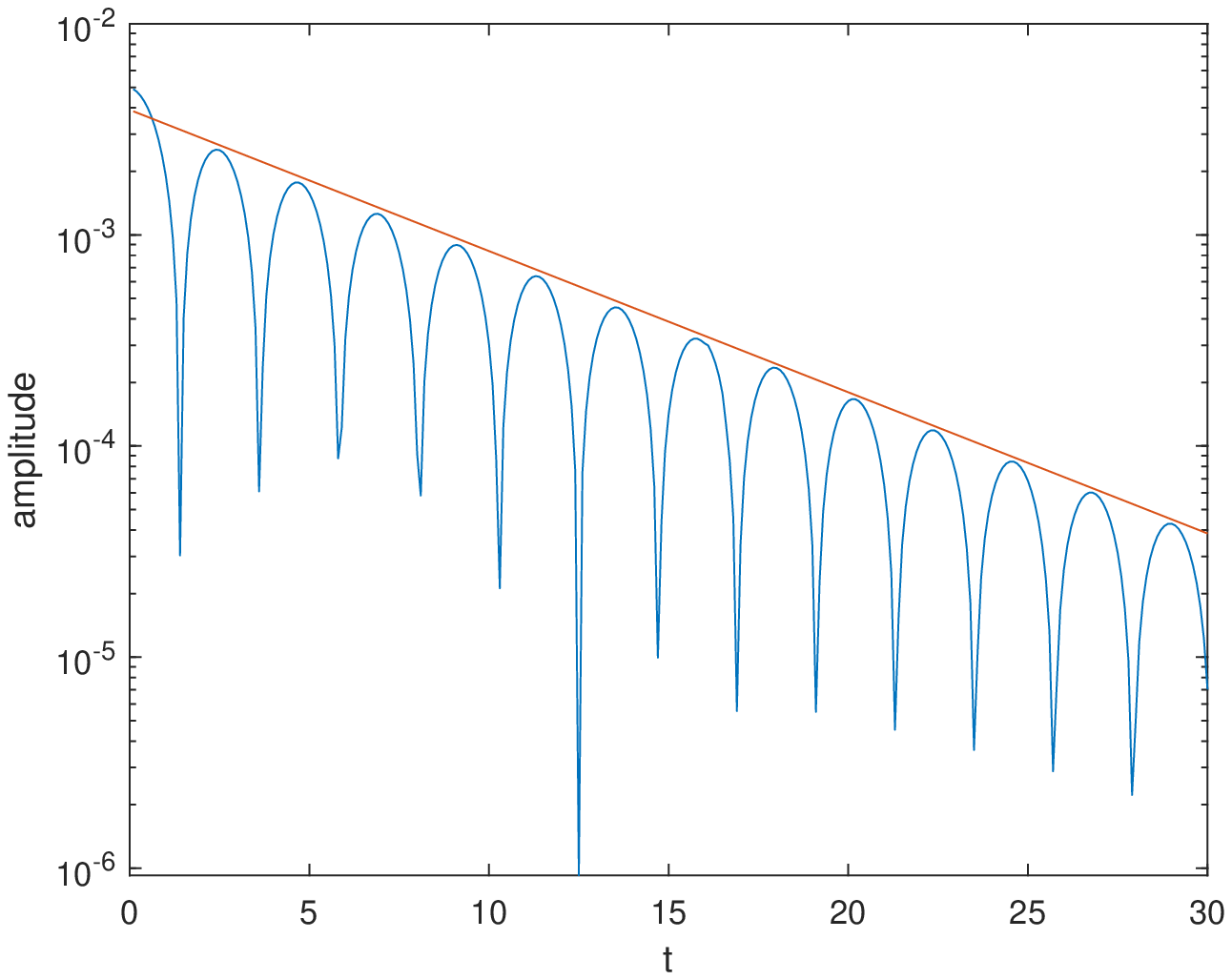}
}
\caption{Electric energy as a function of time for the linear Landau damping. (a) $k=0.3$;\quad (b) $k=0.5$.}
\label{imgL1}
\end{figure}

\begin{figure}[h!]
\centering
\subfigure[]{
\includegraphics[scale=.5]{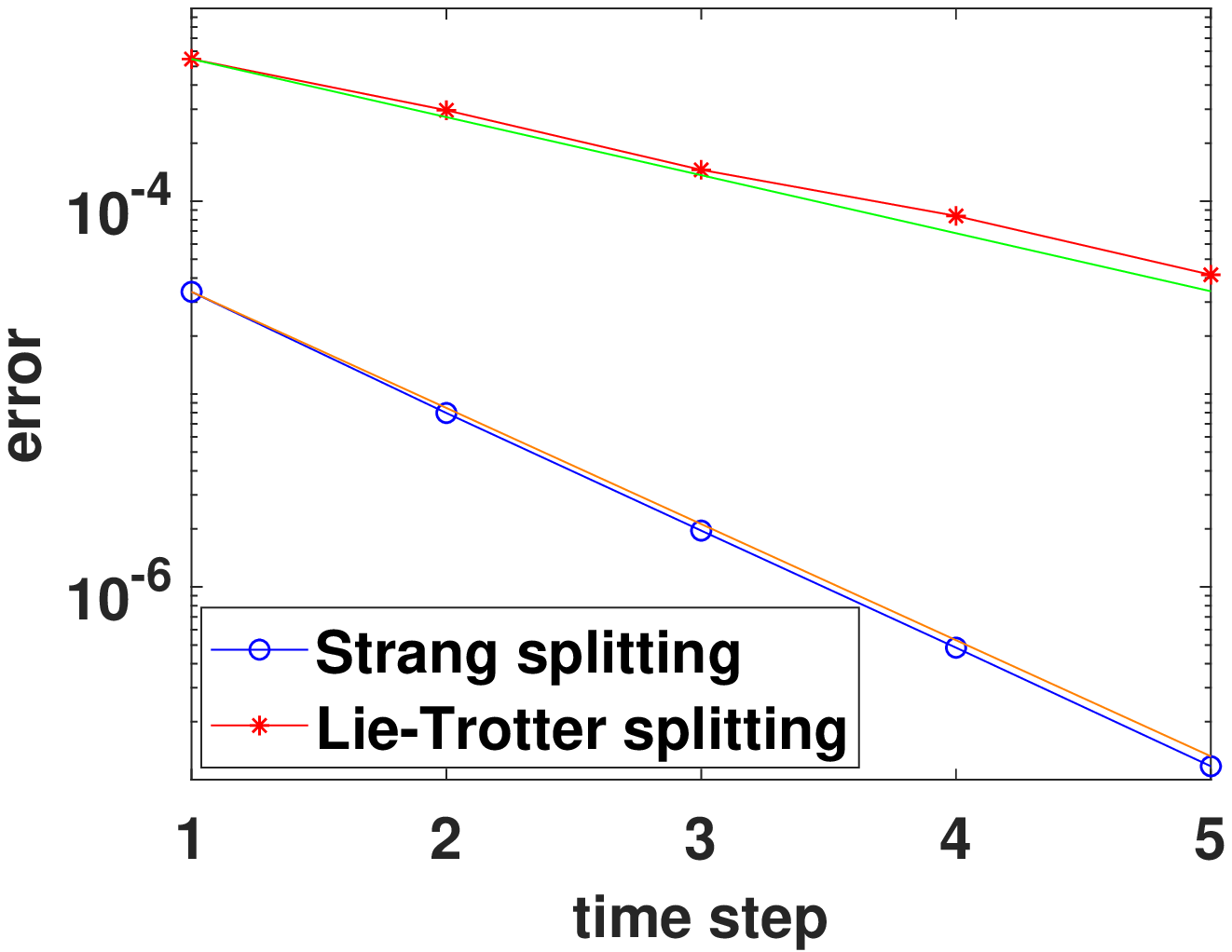}
}
\quad
\subfigure[]{
\includegraphics[scale=.5]{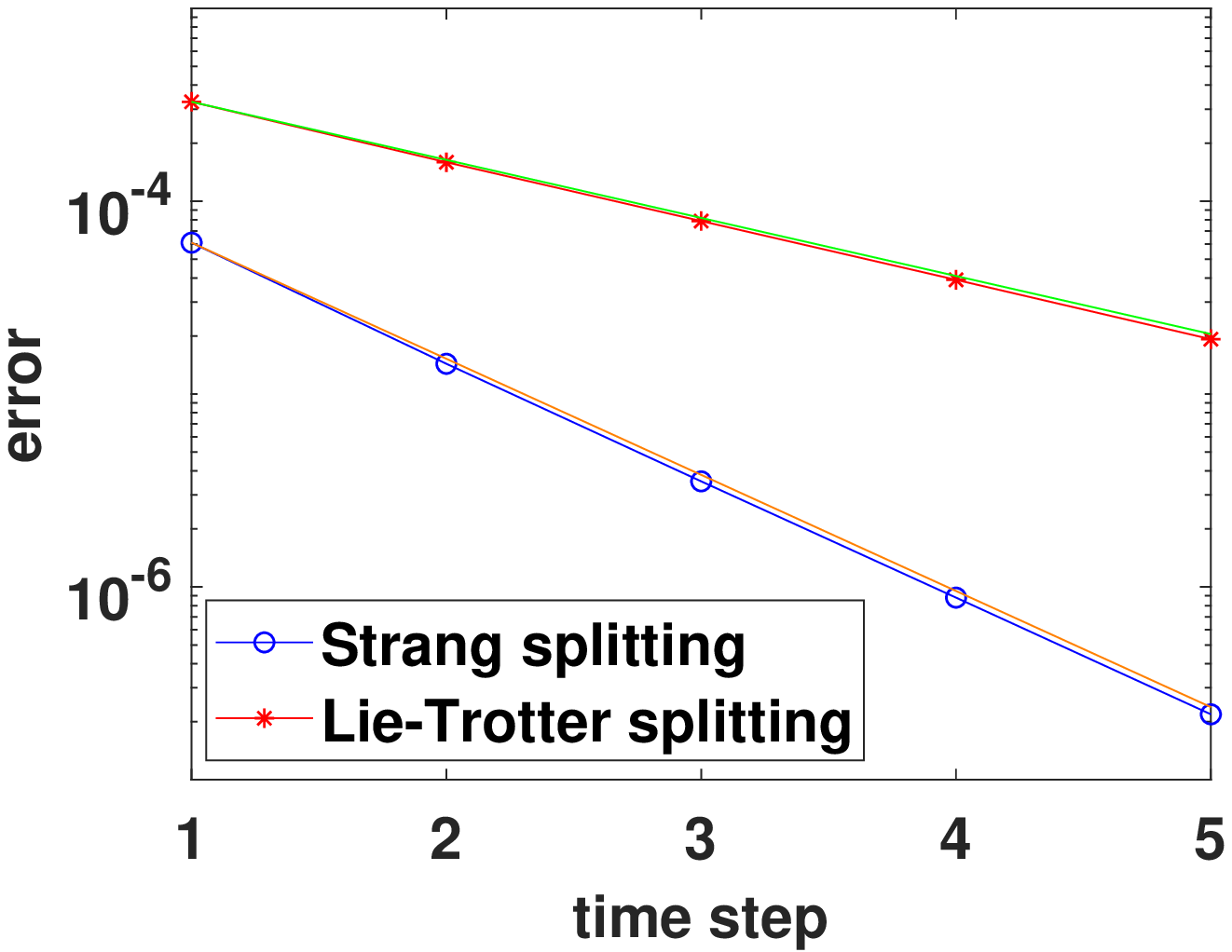}
}
\caption{Convergence rates of numerical solutions with respect to time step. (a): Error of position; (b): Error of velocity. The solid colored lines refer to the corresponding theoretical value.}
\label{imgL2}
\end{figure}

We use the time step $\Delta{t}=0.01$. To solve the Poisson system we use $N_x=128$.
We choose the computation domain as $[0,2\pi/k]\times[-6, 6]$. We are interested in the
evolution of the square root of the electric energy which is measured by ${E}_{d}(t)=\sqrt{\mathbf{\Phi(X)}^{\mathrm T} \mathbb{M} \mathbf{\Phi(X)}}$ with $\mathbf{\Phi(X)}$ and $\mathbb{M}$  mentioned in Section 4.
In Figure \ref{imgL1}, we calculate the value  of $\log({E}_{d}(t))$ along  time $t$  where the red line is the theoretical damping rate. It is observed that the electric energy is exponentially decreasing, and the damping rate coincides with the analytical values which are $\gamma=0.0127$ for $k=0.3$ and $\gamma=0.154$ for $k=0.5$~\citep{Filbet2001,Besse2003,Nicolas2009}.

We check the convergence order of splitting methods in Figure \ref{imgL2} for Landau damping with $k=0.5$. We start from choosing $\Delta{t}_0=0.5$ and once again conduct four runs by decreasing the time step with each run. That is, the time steps can be chosen as $\Delta{t}_i=2^{-(i+1)},i=1,2,3,4$. We use $e_i^x$ to denote the error of position with time step $\Delta{t}_i$. For a fixed final time $T$, the error order can be calculated by $\log(e_j^x/e_{j+1}^x)/\log(\Delta{t}_j/\Delta{t}_{j+1}),j=0,1,2,3,4$. Similarly, we can calculate the error order of velocity. As it shown in the figure, the numerical results by Strang splitting perform much better than the ones by  Lie-Trotter splitting  due to the higher-order accuracy of Strang splitting method. However, the  damping rates by both approaches can match the analytical ones.

{\bf Two-stream instability.} Two-stream is a common instability in plasma physics which occurs when an energetic particle stream injects into a plasma, or a current is set along
the plasma.
The phenomena can be  investigated when there exists  (time-depending) difference of drift velocity between two plasma components. In our simulation, for the initial value we set a perturbation as a vortex creation at the center of concerned domain. That is, the initial datum can be set as
\begin{equation*}
f_0(x,v)=\frac{1}{\sqrt{2\pi}}{v^2}\exp(-\frac{v^2}{2})(1+\alpha\cos(kx)),
\end{equation*}
where perturbation parameter $\alpha=0.01$ and wave number $k=0.5$. The computation domain is $[0,2\pi/k]\times[-5, 5]$.  We take $\Delta{t}$ and $N_x$ as shown in the above experiments.

\begin{figure}[h!]
\centering
\subfigure{
\includegraphics[scale=.5]{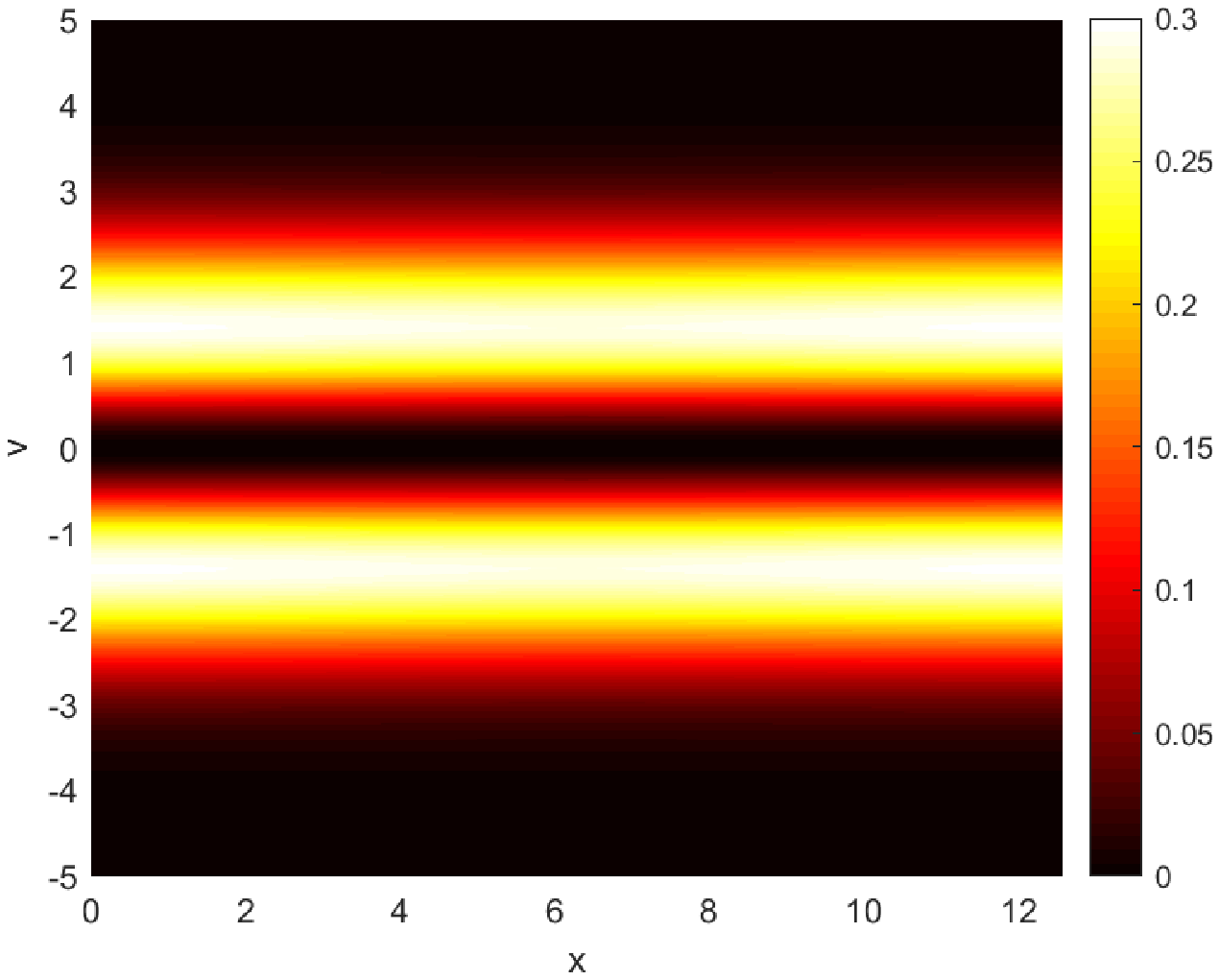}
}
\quad
\subfigure{
\includegraphics[scale=.5]{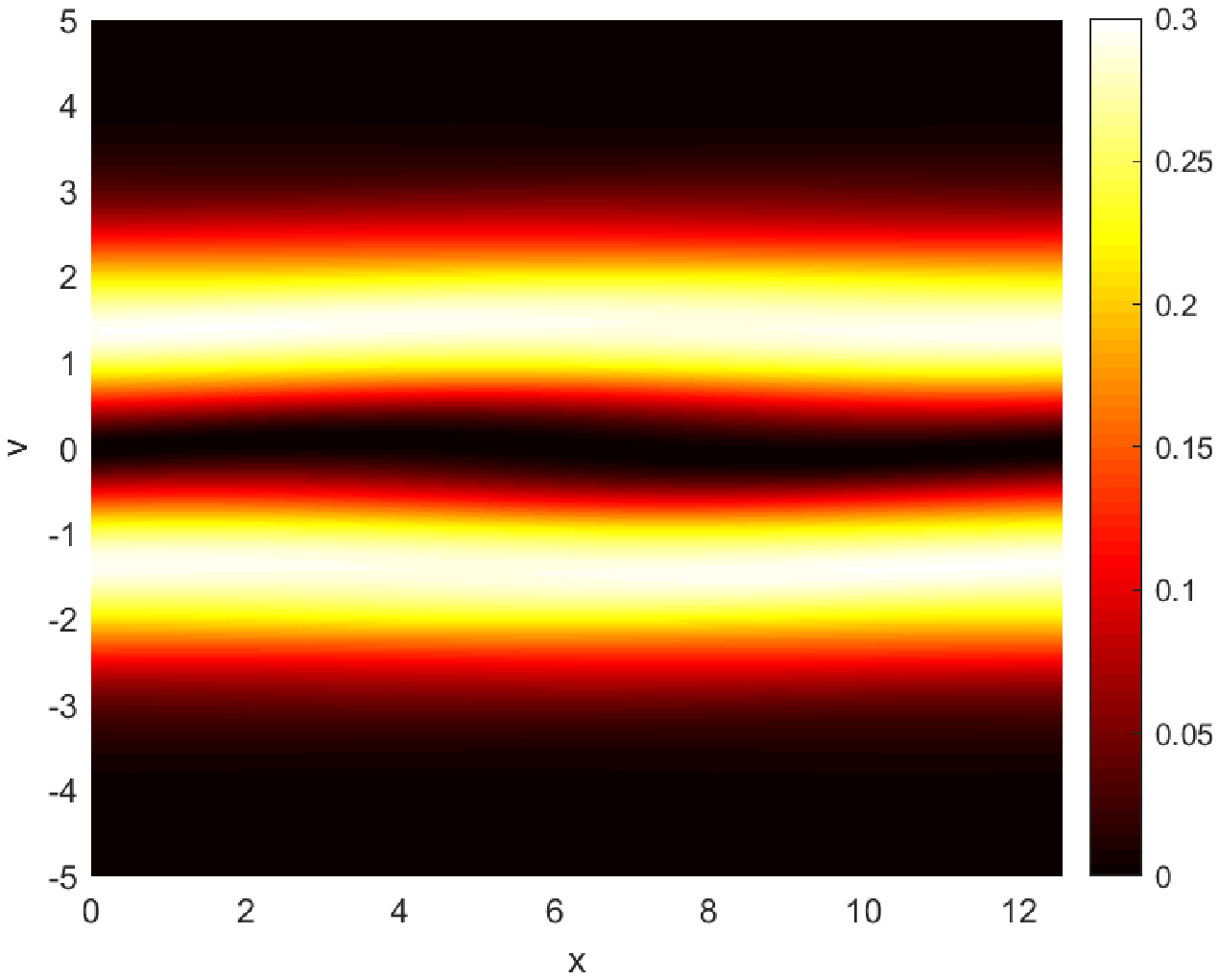}
}
\quad
\subfigure{
\includegraphics[scale=.5]{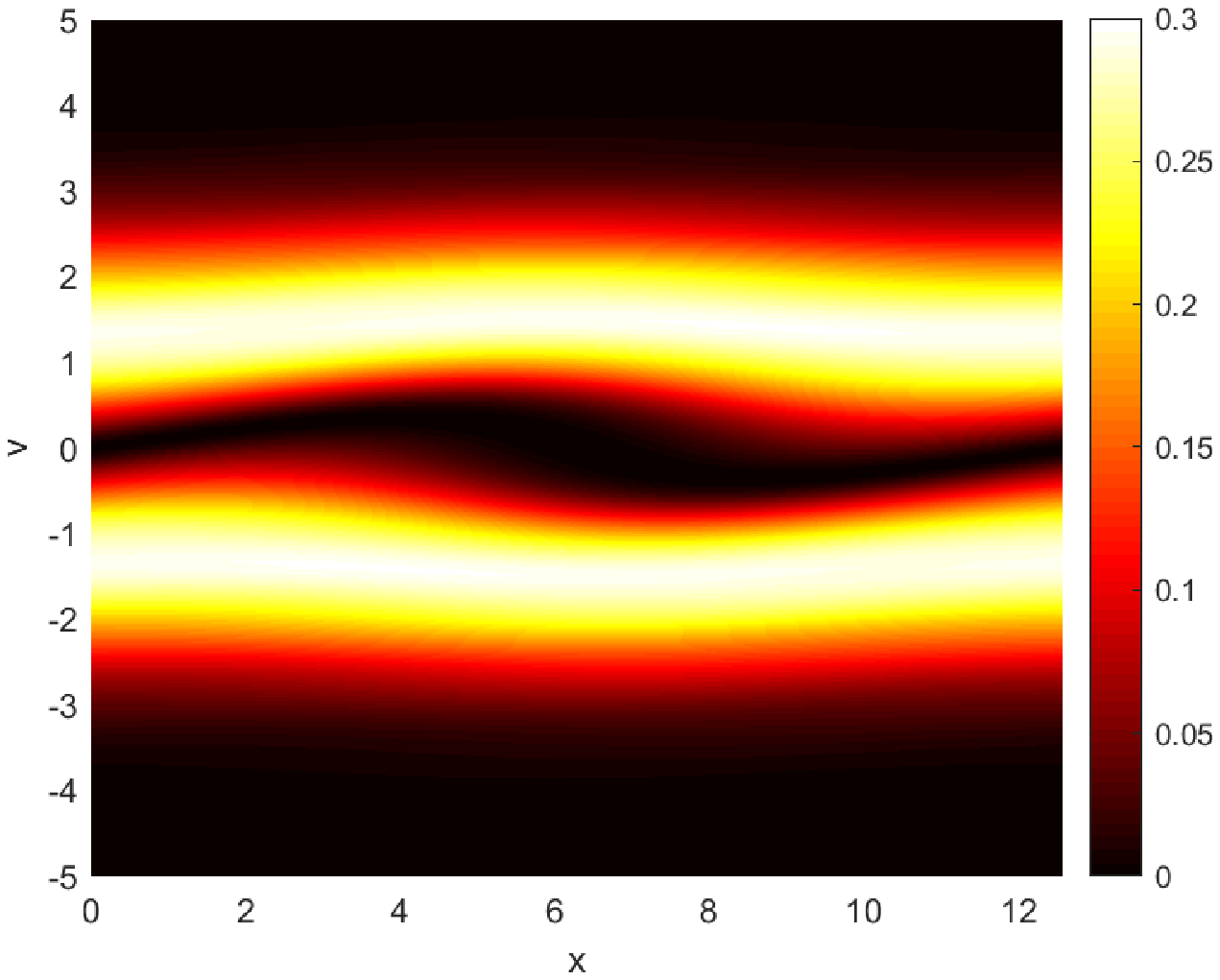}
}
\quad
\subfigure{
\includegraphics[scale=.5]{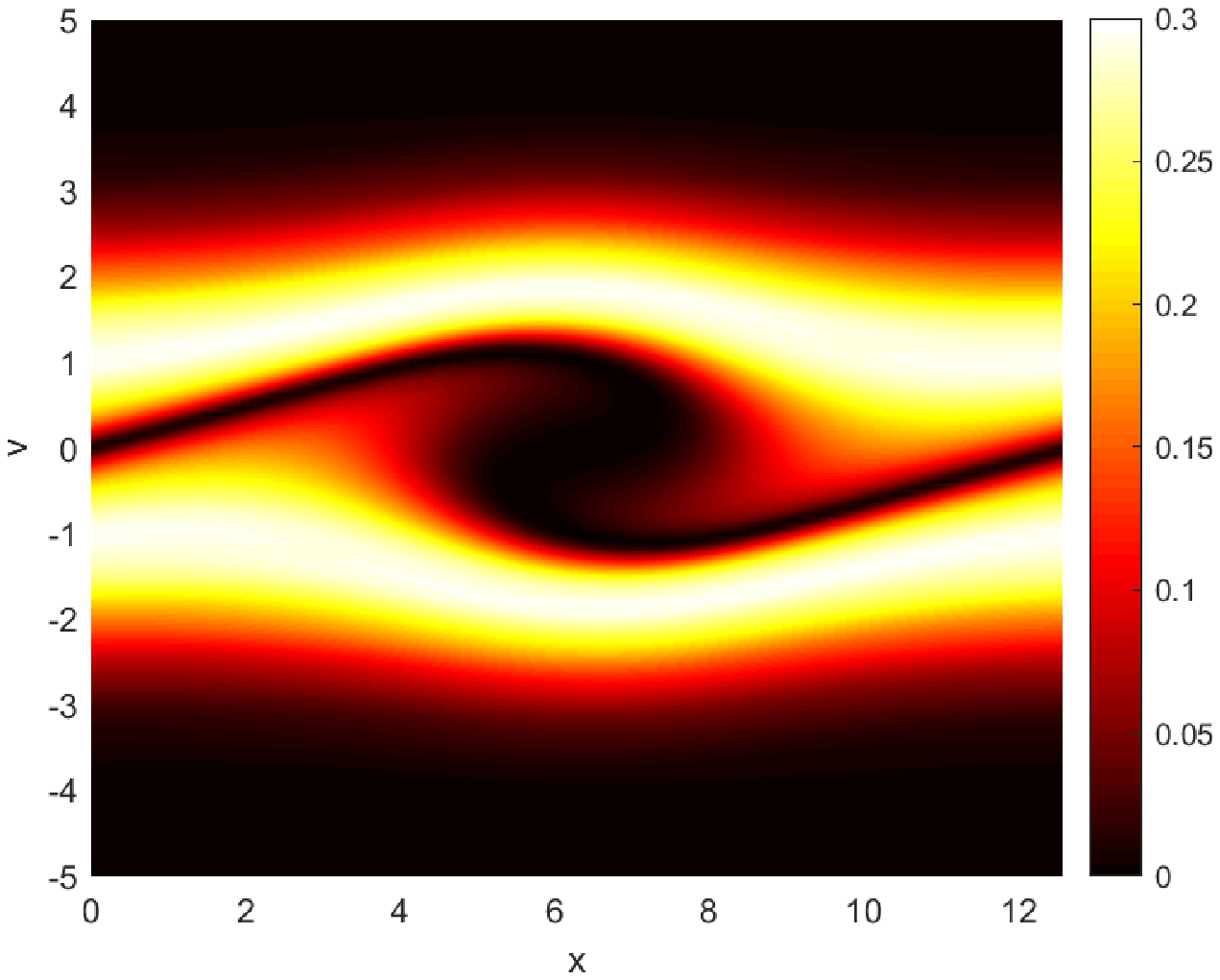}
}
\caption{Evolution of the distribution $f$ for time $t=0,10,15,20$.}
\label{imgT1}
\end{figure}

\begin{figure}[h!]
\centering
\subfigure[]{
\includegraphics[scale=.4]{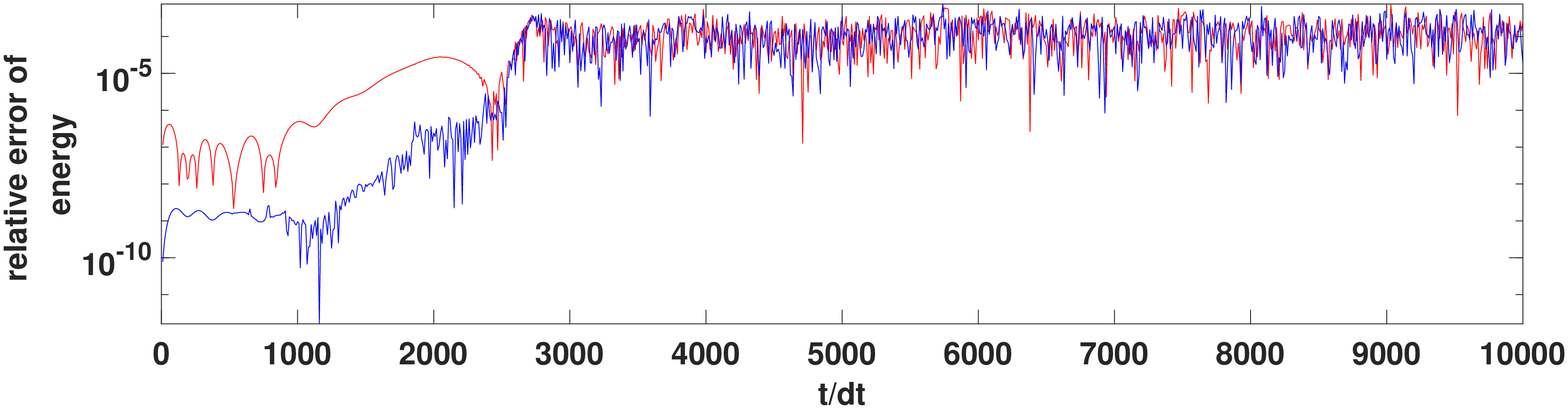}
}
\quad
\subfigure[]{
\includegraphics[scale=.4]{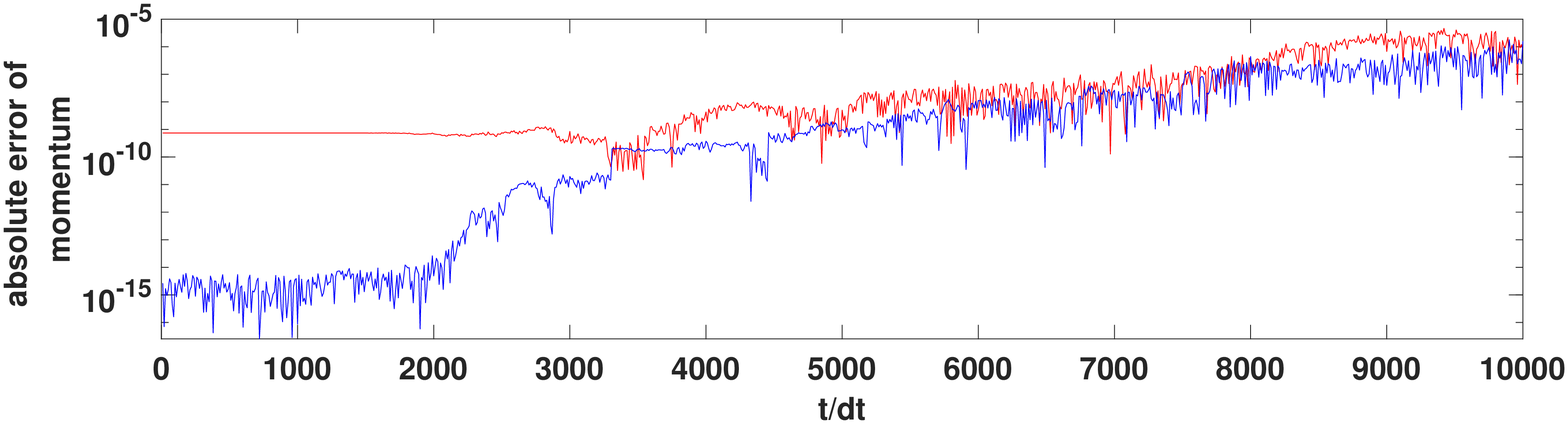}
}
\caption{Time evolutions of conservative quantities  for Two-stream instability. (a): Error of energy; (b): Error of  momentum.  The red line shows the error by the Lie-Trotter splitting while the blue line refers to the ones by the Strang splitting.}
\label{imgT2}
\end{figure}

In Figure \ref{imgT1}, we plot the distribution function $f$ for time $t=0,10,15,20$.
It is observed that the instability appears quickly, and the obvious vortex structure can be captured clearly.
This coincides with the instability phenomenon shown in literature~\citep{Besse2003,Nicolas2009}.

In Figure \ref{imgT2}, we plot the total energy and momentum. We can observe that our method can preserve the energy and momentum well over  long-time even when the instability occurs.

{\bf{Bump-on-tail instability.}} Bump-on-tail instability is a fundamental example of wave-particle interaction, the instability  can be investigated when the wave and particles interact.
In this test, for studying the Bump-on-tail instability we take the following initial
conditions~\citep{Shoucri2011}
\begin{equation*}
f_0(x,v)=(1+\alpha\cos(kx))d_0(v),
\end{equation*}
where perturbation parameter $\alpha=0.04$ and wave number $k=0.3$. In the above initial value, $d_0$ denotes a bump  on the tail of  distribution function  formed by energetic particles,  which has the following expression
\begin{equation*}
d_0=\frac{1}{\sqrt{2\pi}}\left(0.9\exp\left(-\frac{v^2}{2}\right)+0.2\exp(-2(v-4.5)^2)\right).
\end{equation*}

\begin{figure}[h!]
\centering
\subfigure{
\includegraphics[scale=.5]{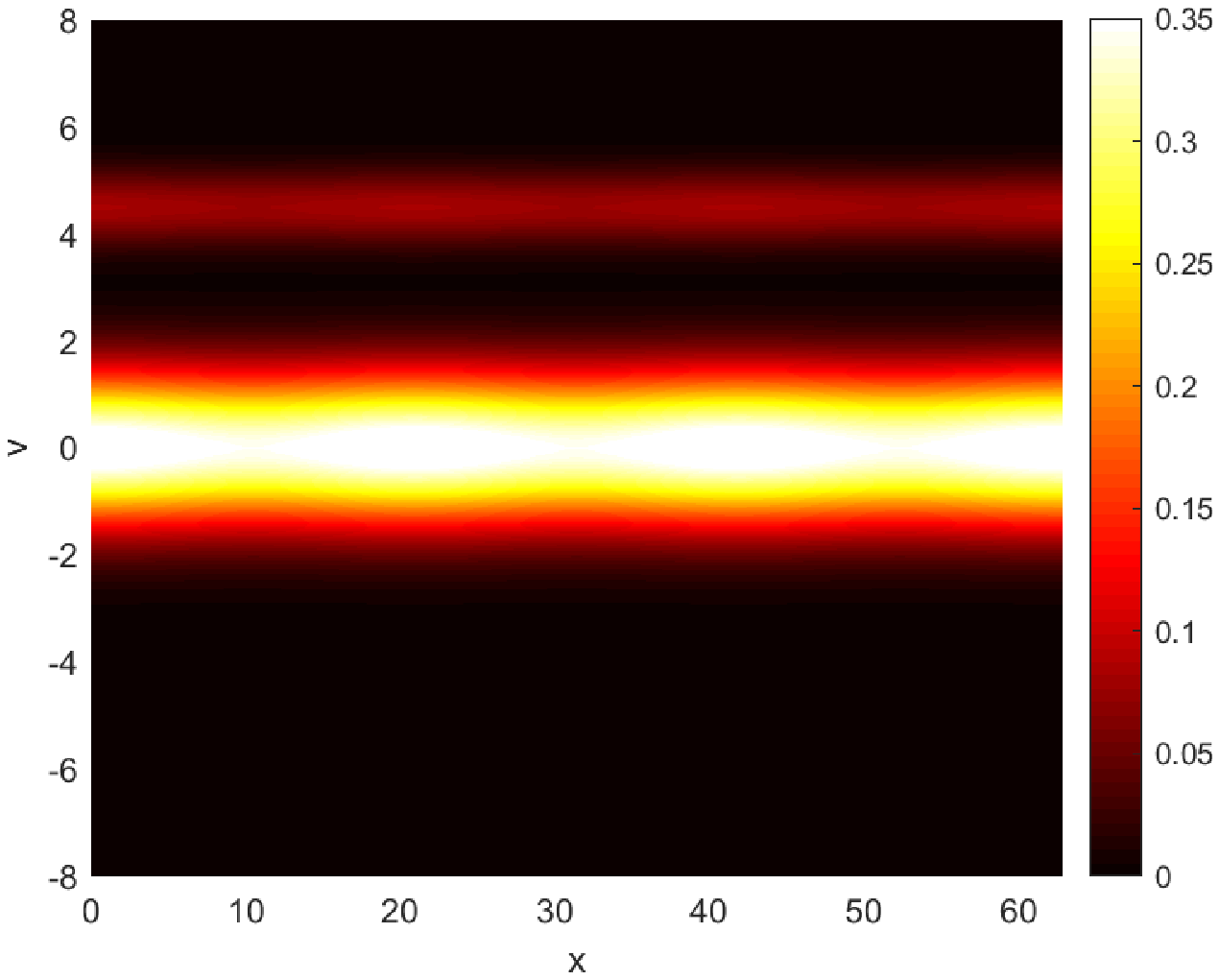}
}
\quad
\subfigure{
\includegraphics[scale=.5]{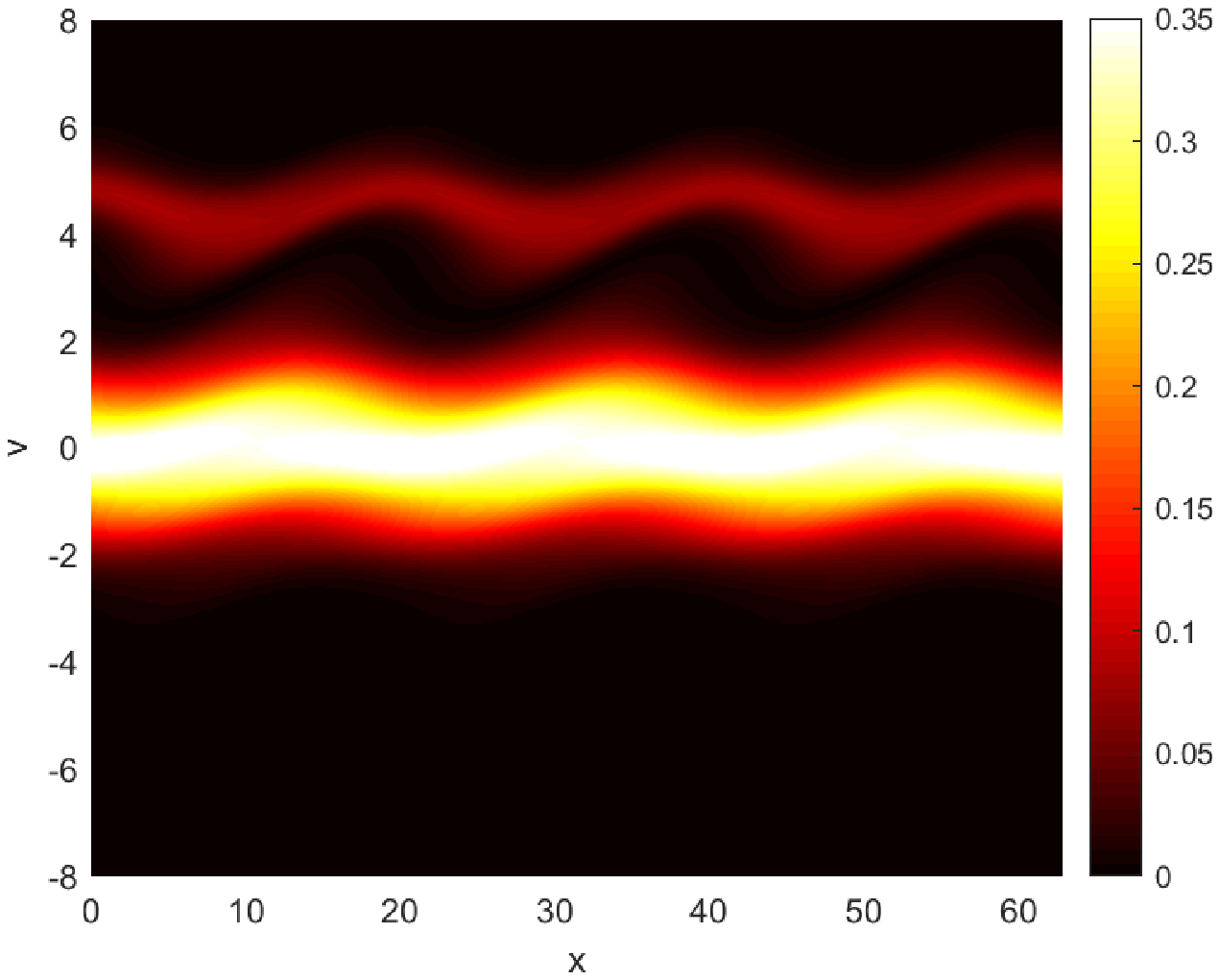}
}
\quad
\subfigure{
\includegraphics[scale=.5]{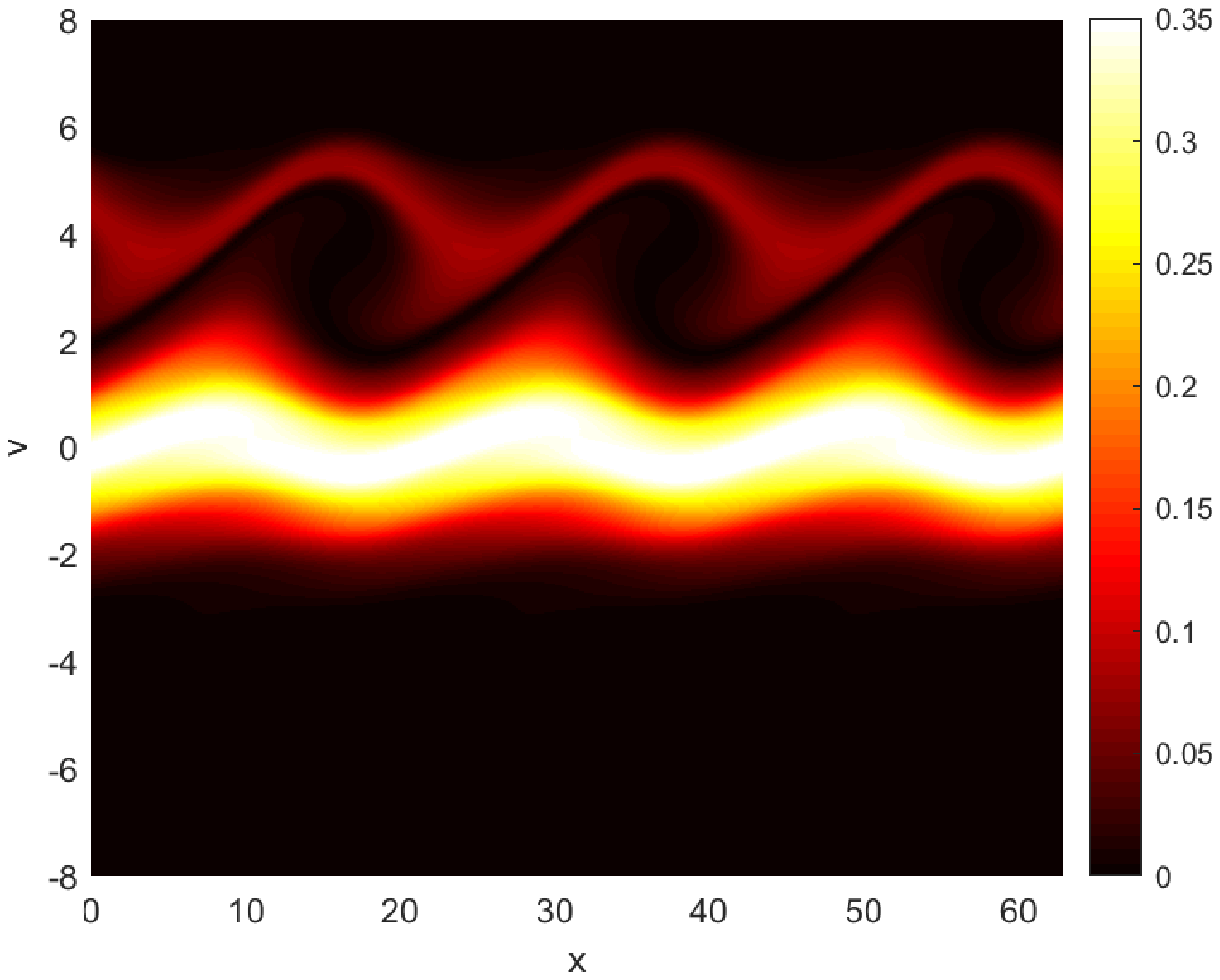}
}
\quad
\subfigure{
\includegraphics[scale=.5]{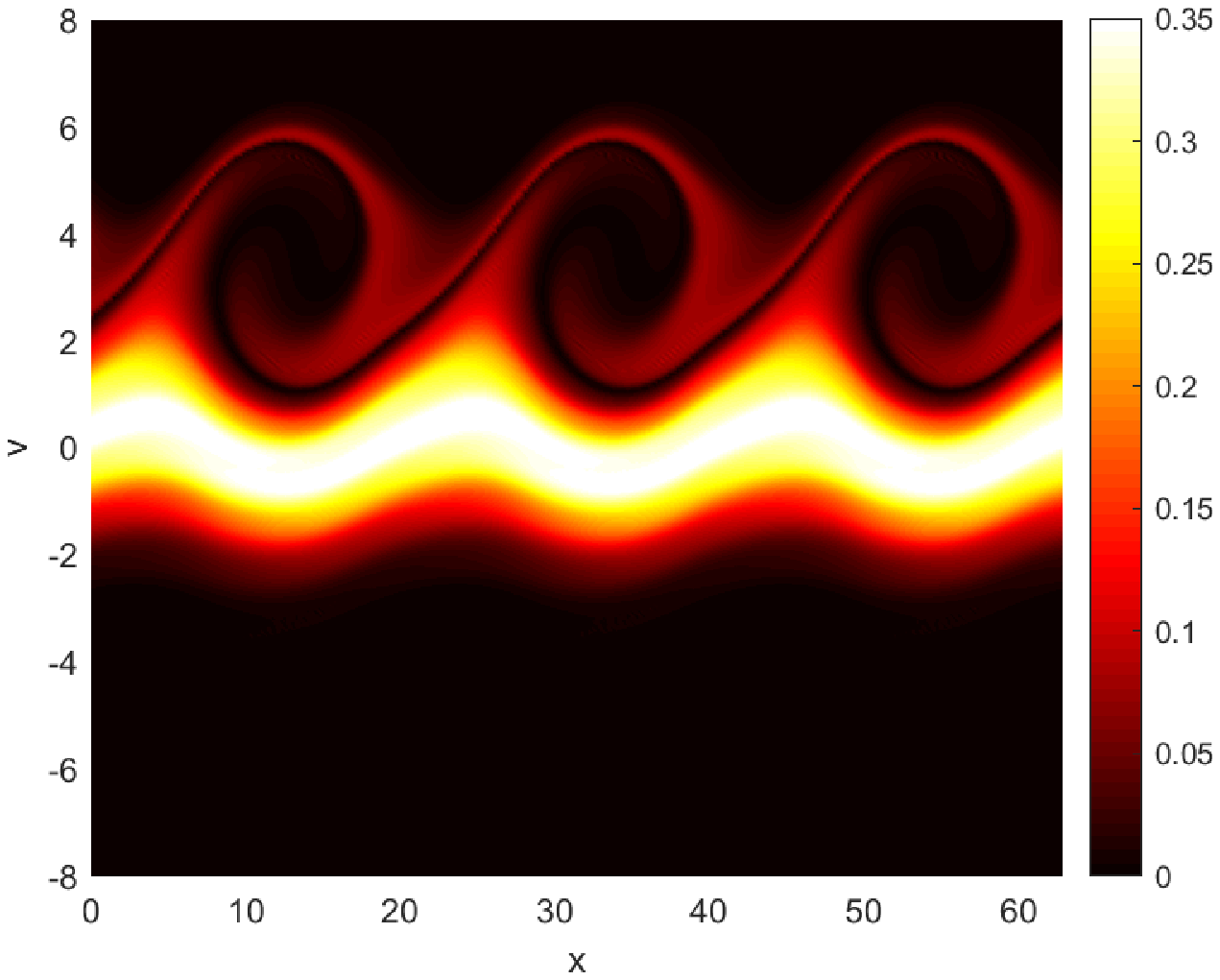}
}
\caption{Time evolution of the distribution $f$ for time $t=0,10,15,20$.}
\label{imgB1}
\end{figure}

In this experiment, computation domain $[0,6\pi/k]\times[-8, 8]$ is considered. We choose $\Delta{t}=0.01$ and $N_x=256$ to simulate this problem.
In Figure \ref{imgB1}, we plot the distribution function $f$ for time $t=0,10,15,20$.  It is observed that our numerical discretization can describe the  evolution of instability and the occurrence phenomena. Clearly, it is investigated that there appear three vortices which extending arms  embrace the neighbouring vortices.

\begin{figure}[h!]
\centering
\includegraphics[scale=.5]{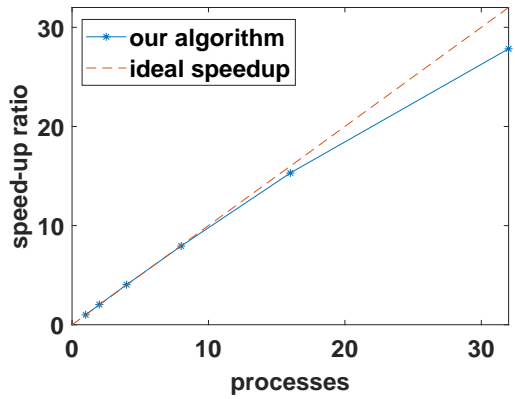}
\caption{Speed-up of numerical computing for simulating Bump-on-tail instability on LSSC-IV. The computation is done by taking $10^5$ particles in phase space.}
\label{imgSR}
\end{figure}

It should  be noticed that implementing our numerical method involves computations in a large amount of cycles which can be treated in parallel. As follows, we show the efficiency of our numerical computation together with  the parallel approach. Based on MPI, we take two natural parallel strategies for the collective communications in the cycle. For particle deposition part, the computational domain is subdivided into regions while for pushing particle part, the particles are divided into groups. For the first strategy, we can also improve it by using, for instance, the adaptive scheme presented in~\cite{Mehrenberger2006}.
In Figure \ref{imgSR}, we present the speed-up comparison for our numerical computation and the ideal ones. By comparison, it can be seen  that our speed-up result is quite good due to  the so-called cache effects.

\subsection{Test problems in 2+2-dimensional phase space}
In this experiment, we consider the Vlasov--Poisson system in 2+2-dimension which means 2D in space and also 2D in velocity one. We assign $N_x=N_y=256$ for the initial particle grid and $250$ particles per cell. This  means that we take over $10^7$ particles for the simulation.

{\bf{Diocotron instability.}} The diocotron instability is a plasma instability which occurs when two sheets of charges slipping past each other. This stability is usually described by the guiding center model~\citep{Levy1965,Filbet2016,Ameres2018} and has been extensively simulated~\citep{Nicolas2014,Filbet2016,Filbet2018,Ameres2018,Chartier2020}. Here, we consider the model of two dimensional Vlasov--Poisson system together with an external magnetic field. This system can be taken to model the instability which is usually encountered in magnetic fusion devices such as tokamaks. In this example, we also consider the effect of magnetic field to the instability. Thus, we consider the case with strong magnetic field. As the instability usually is explained as the analog of the Kelvin-Helmholtz instability in fluid mechanics, it can be observed directly in the nature~\citep{Peurrung1993}.  In this situation,  the energy of system is dissipated as the two surface waves propagate in two opposite directions with one flowing over the other. By investigation,  there appears the vortex structure in the surface of distribution function  when the instability happens~\citep{Filbet2016,Filbet2018,Ameres2018}.

The Vlasov equation in this simulation is with an external magnetic field $\mathbf{B}_{\rm{ext}}$ which reads
\begin{equation*}
\varepsilon\frac{\partial f}{\partial t}+\mathbf{v}\cdot\frac{\partial f}{\partial\mathbf{x}}+{(\mathbf{E}(t,\mathbf{x})+\frac{1}{\varepsilon}\mathbf{v}\times\mathbf{B}_{\rm{ext}}(t,\mathbf{x}))}\cdot\frac{\partial f}{\partial\mathbf{v}}=0.
\end{equation*}
Here, $\varepsilon$ is the strength of magnetic field. As the parameter $\varepsilon$ is usually small, the existence of term $\varepsilon$ in front of the time derivative of $f$ requires the numerical simulations  to have long-time stability.
There have been many numerical methods for studying this instability. Among them,  the asymptotic preserving schemes are introduced in~\citep{Filbet2018}, and the uniformly accurate methods are introduced in~\citep{Chartier2020}.

We take the initial distribution function as
\begin{equation*}
f_0(\mathbf{x},\mathbf{v})=\frac{d_0(\mathbf{x})}{2\pi}\exp(-\frac{{\left\|\mathbf{v}\right\|}^2}{2}),\quad \mathbf{x}=(x,y)\in\mathbb{R}^2,
\end{equation*}
where the initial density is
\begin{equation*}
d_0(\mathbf{x})=
\begin{cases}
(1+{\alpha}\cos(l\theta))\exp(-4(\left\|\mathbf{x}\right\|-6.5)^2)&\mbox{if $r^-\leq\left\|\mathbf{x}\right\|\leq{r}^+$,} \\
0&\mbox{otherwise,}
\end{cases}
\end{equation*}
with  $\theta={\rm{atan}}(y/x)$ and $l$ the number of vortices. In our simulation, we take $r^-=5,r^+=8,\alpha=0.2$ and $\mathbf{B}_{\rm{ext}}=(0,0,1)$.

\begin{figure}[h!]
\centering
\subfigure{
\includegraphics[scale=.4]{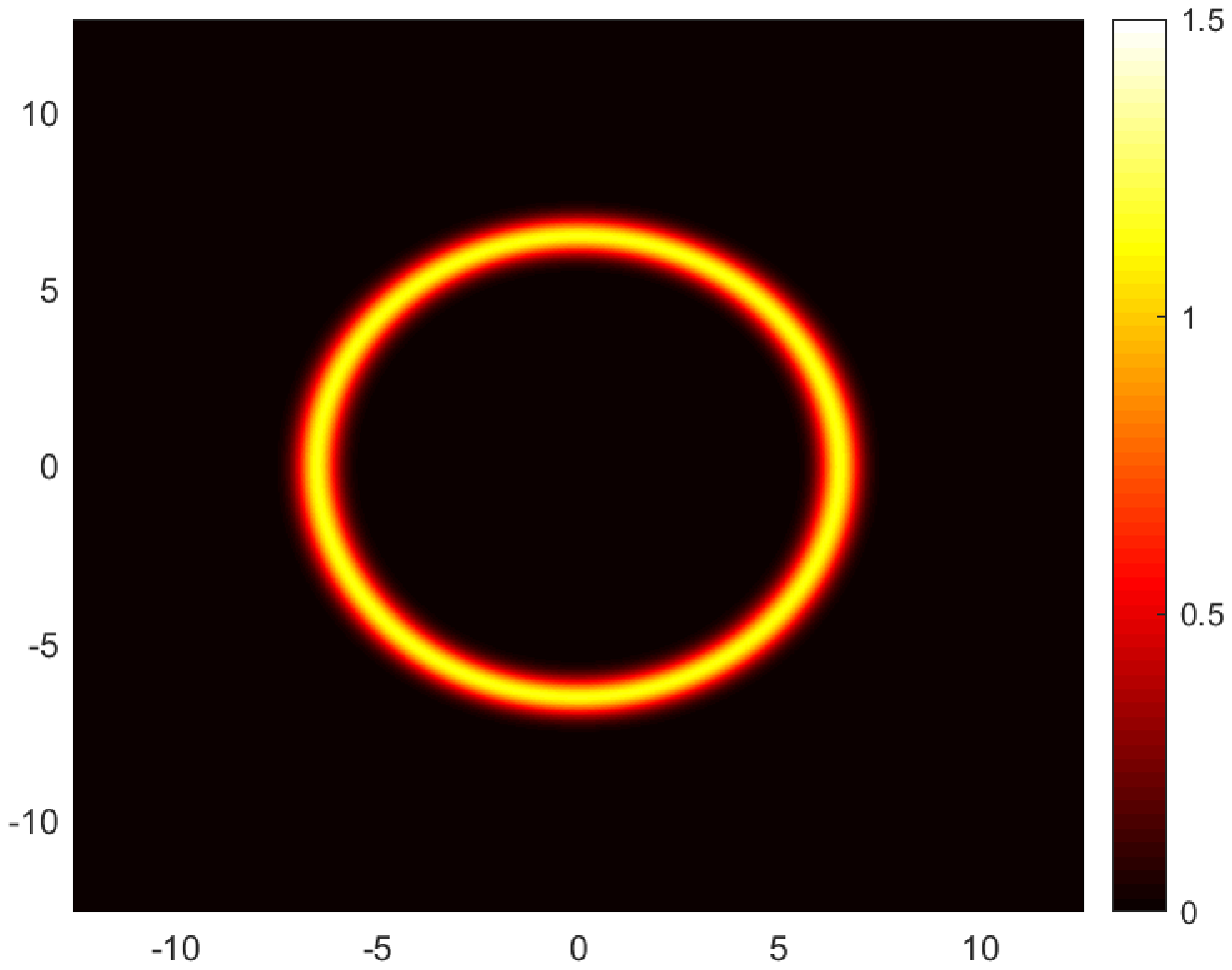}
}
\quad
\subfigure{
\includegraphics[scale=.4]{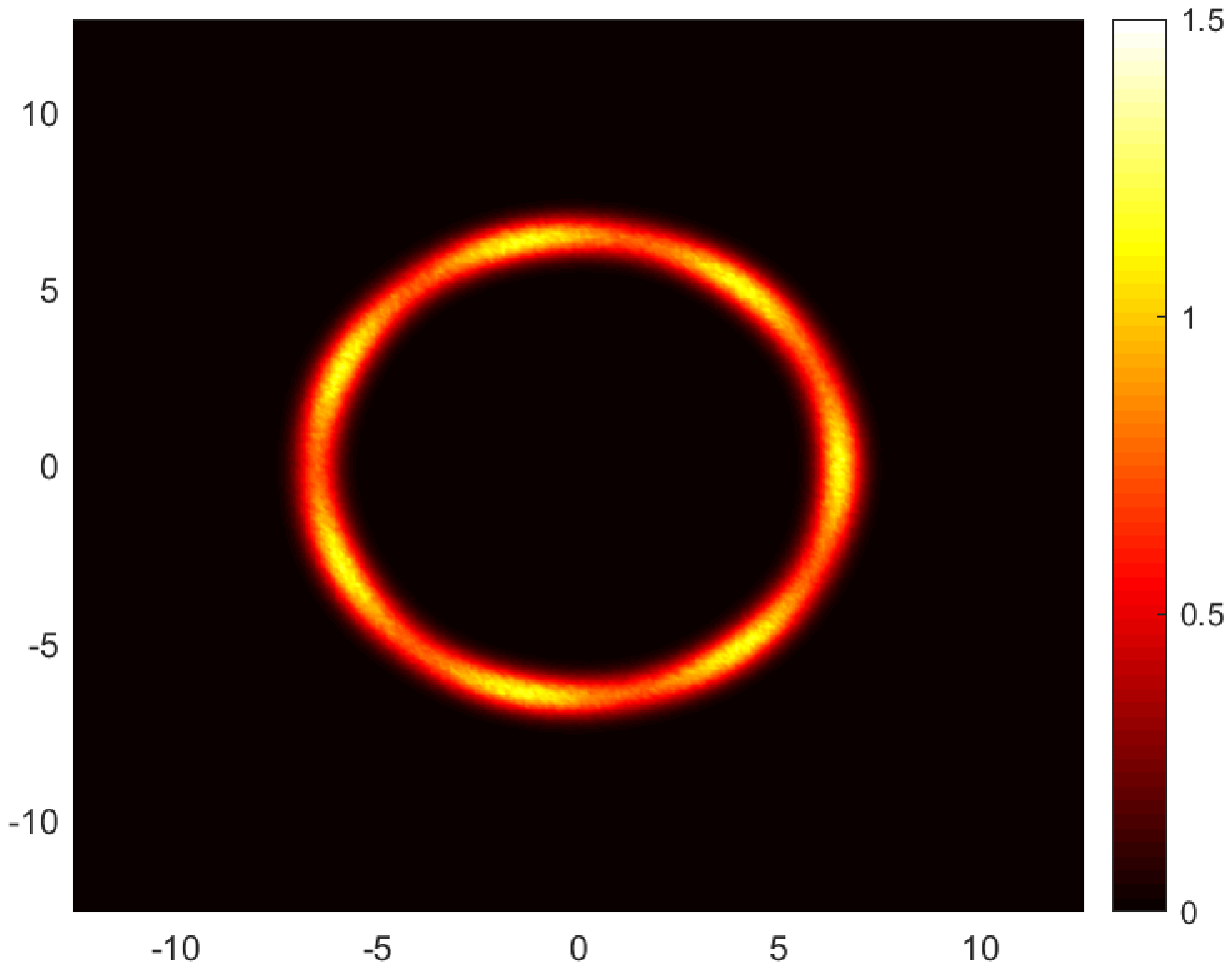}
}
\quad
\subfigure{
\includegraphics[scale=.4]{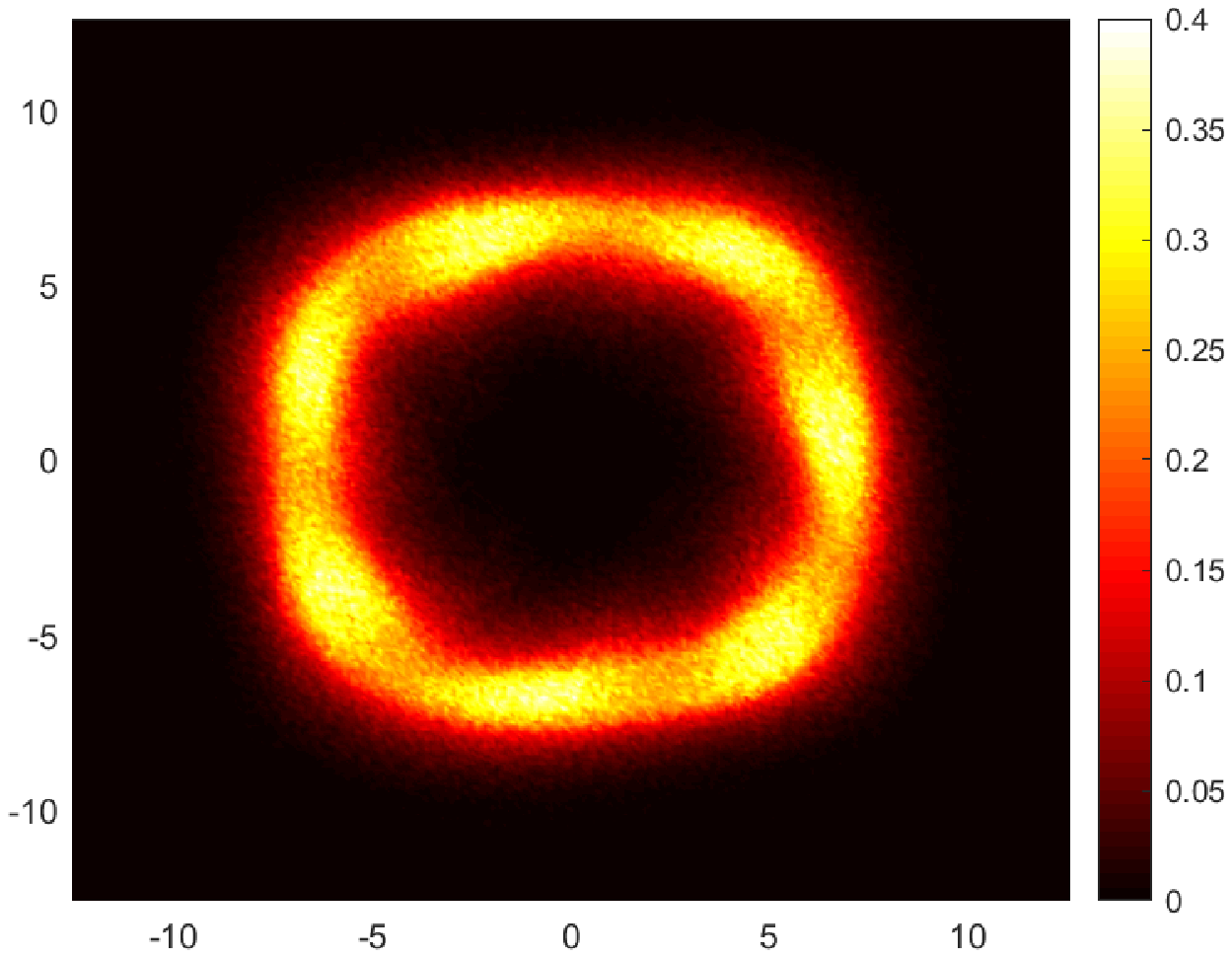}
}
\quad
\subfigure{
\includegraphics[scale=.4]{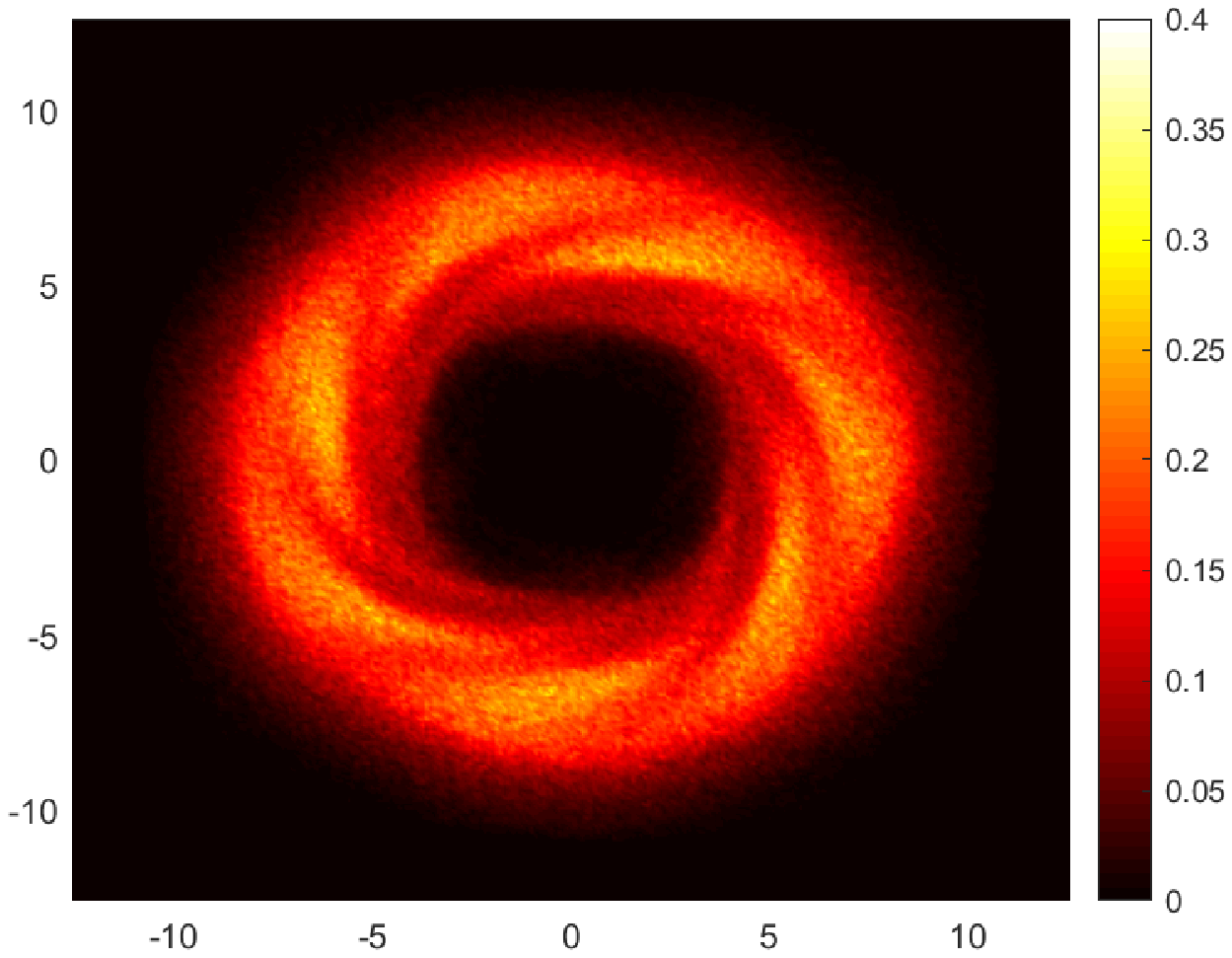}
}
\caption{Time evolution of the density $\rho$ for time $t=0,0.1,15,30$ and $l=7$ with  strength of magnetic field $\varepsilon=1$.}
\label{imgD1}
\end{figure}

In Figure \ref{imgD1}, we consider the case where $\varepsilon$ is taken as $\varepsilon=1$. We use $\Delta{t}=0.1$ and $l=7$ for this case. It can be observed that the plasma is not well confined. This can be explained that  the  enforced magnetic field is not strong enough. This phenomena  is also performed in~\citep{Filbet2016}. Though the vortices will not occur, it still can be observed that for the case of $l=7$ there are seven clear clusters changing along the time.

\begin{figure}[h!]
\centering
\subfigure[]{
\includegraphics[scale=.4]{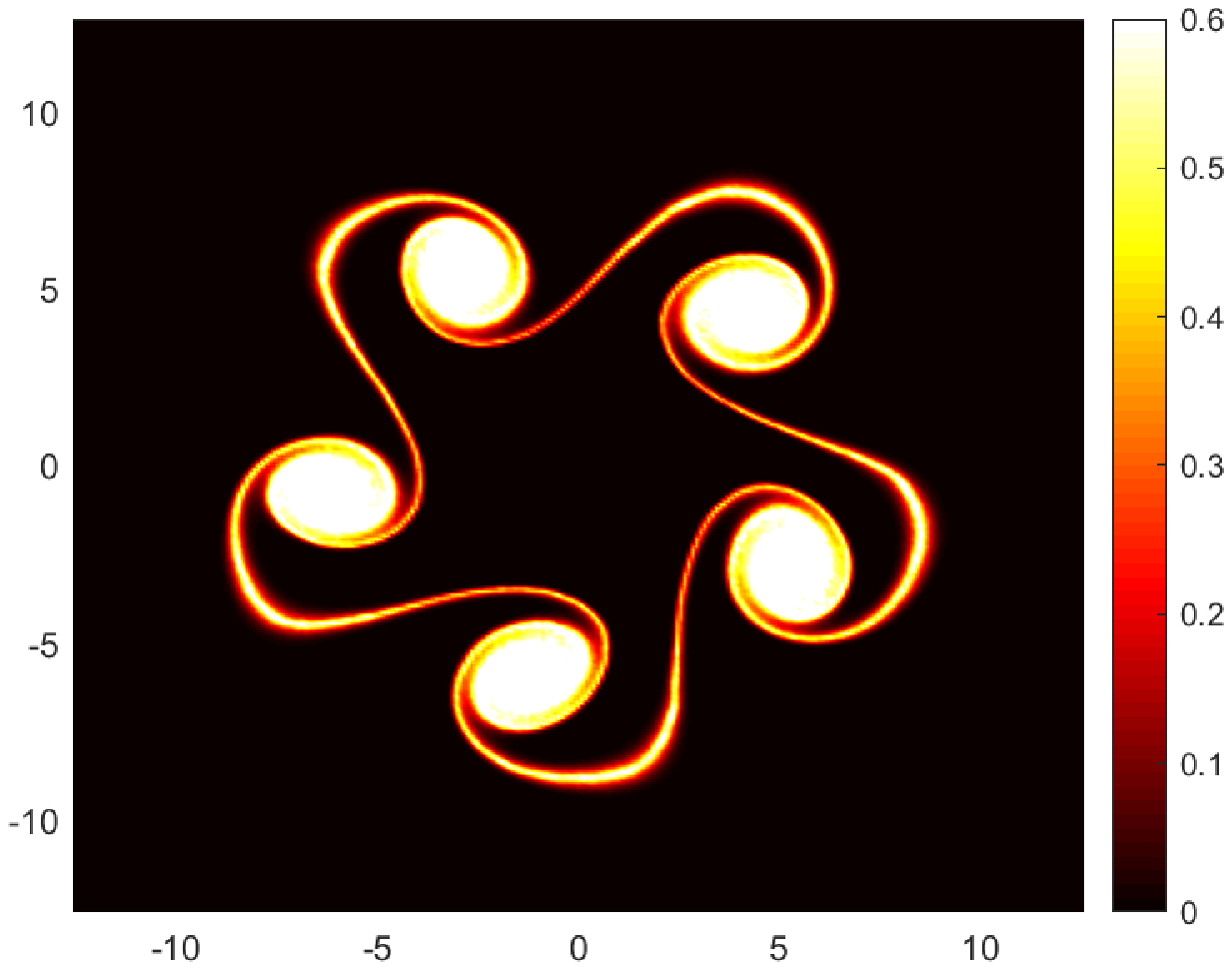}
}
\quad
\subfigure[]{
\includegraphics[scale=.4]{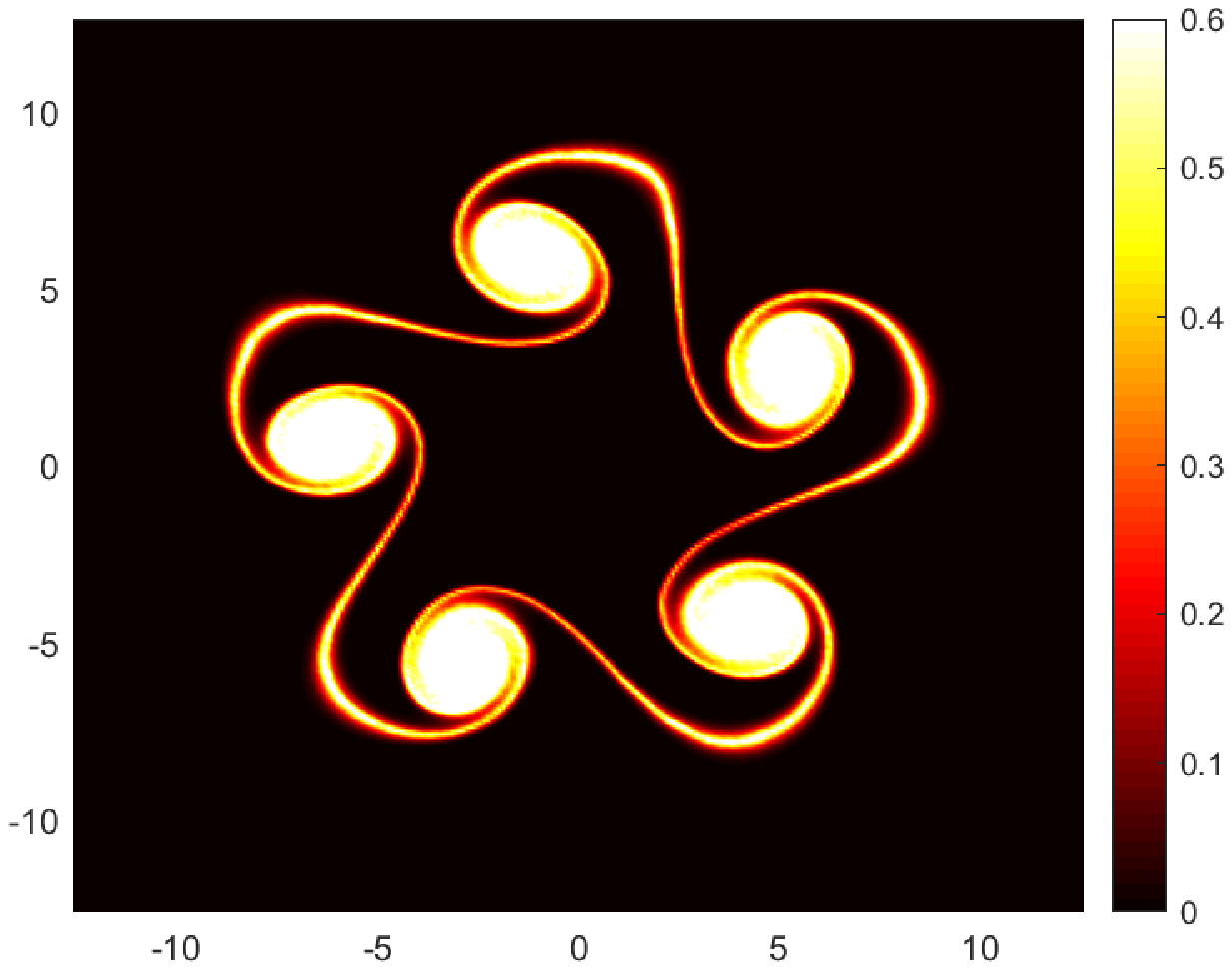}
}
\quad
\subfigure[]{
\includegraphics[scale=.4]{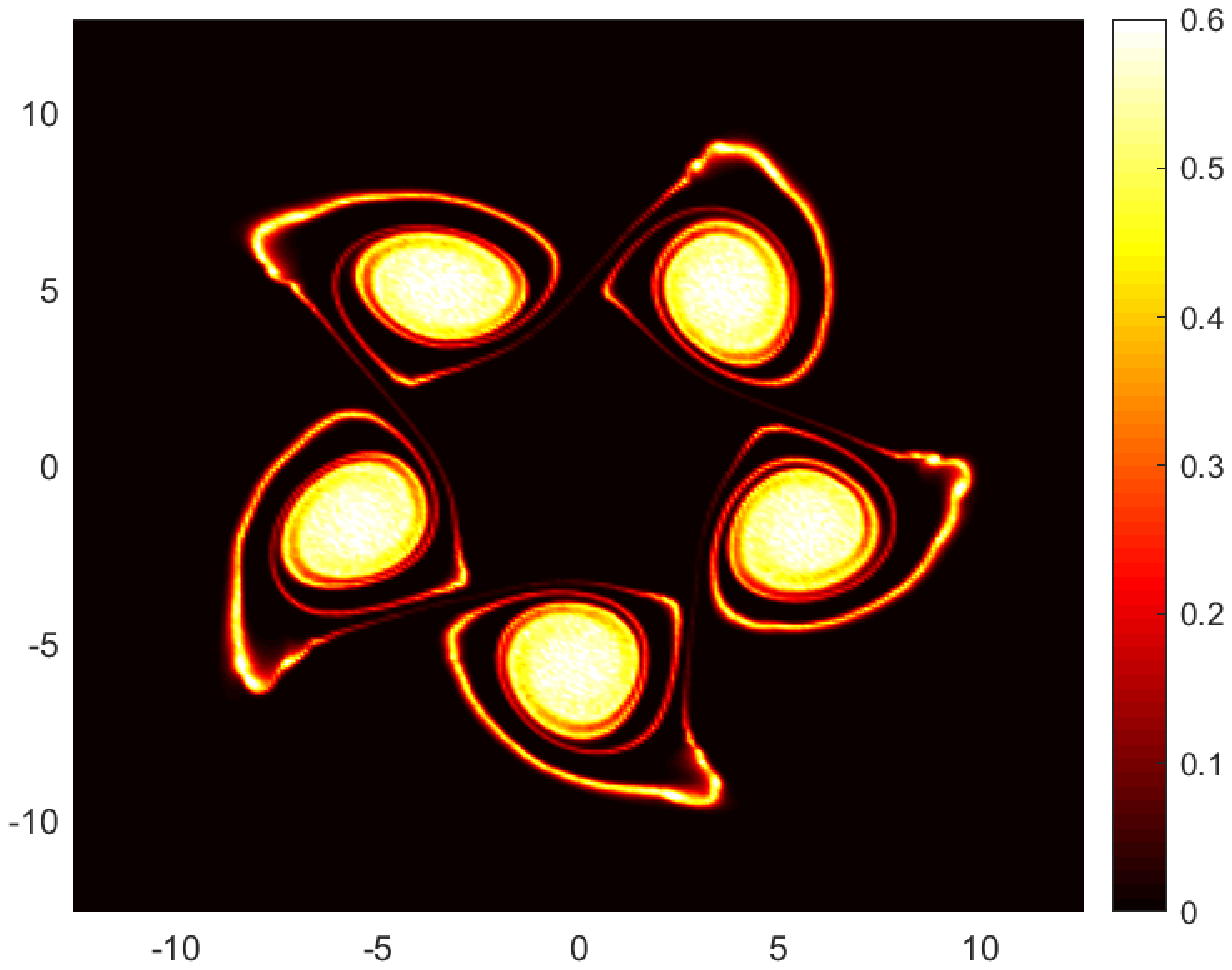}
}
\quad
\subfigure[]{
\includegraphics[scale=.4]{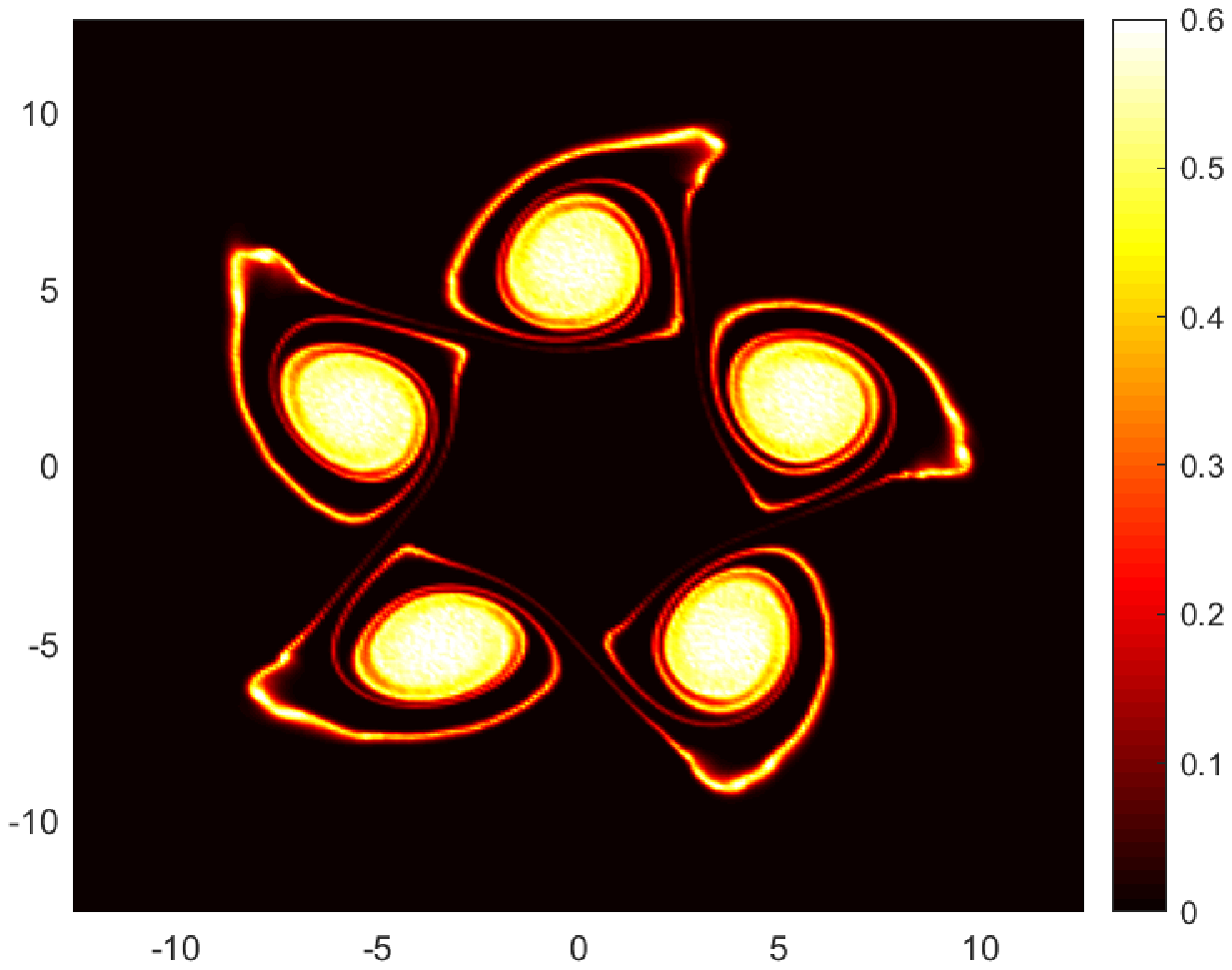}
}
\caption{Time evolution of the density $\rho$ along the time. (a): $\mathbf{B}_{\rm{ext}}=(0,0,1)$, $t=30$; (b): $\mathbf{B}_{\rm{ext}}=(0,0,-1)$, $t=30$; (c): $\mathbf{B}_{\rm{ext}}=(0,0,1)$, $t=50$; (d) $\mathbf{B}_{\rm{ext}}=(0,0,-1)$, $t=50$. Here $l=5$ and the strength of magnetic field  is $\varepsilon=0.1$.}
\label{imgD2}
\end{figure}

\begin{figure}[h!]
\centering
\subfigure{
\includegraphics[scale=.4]{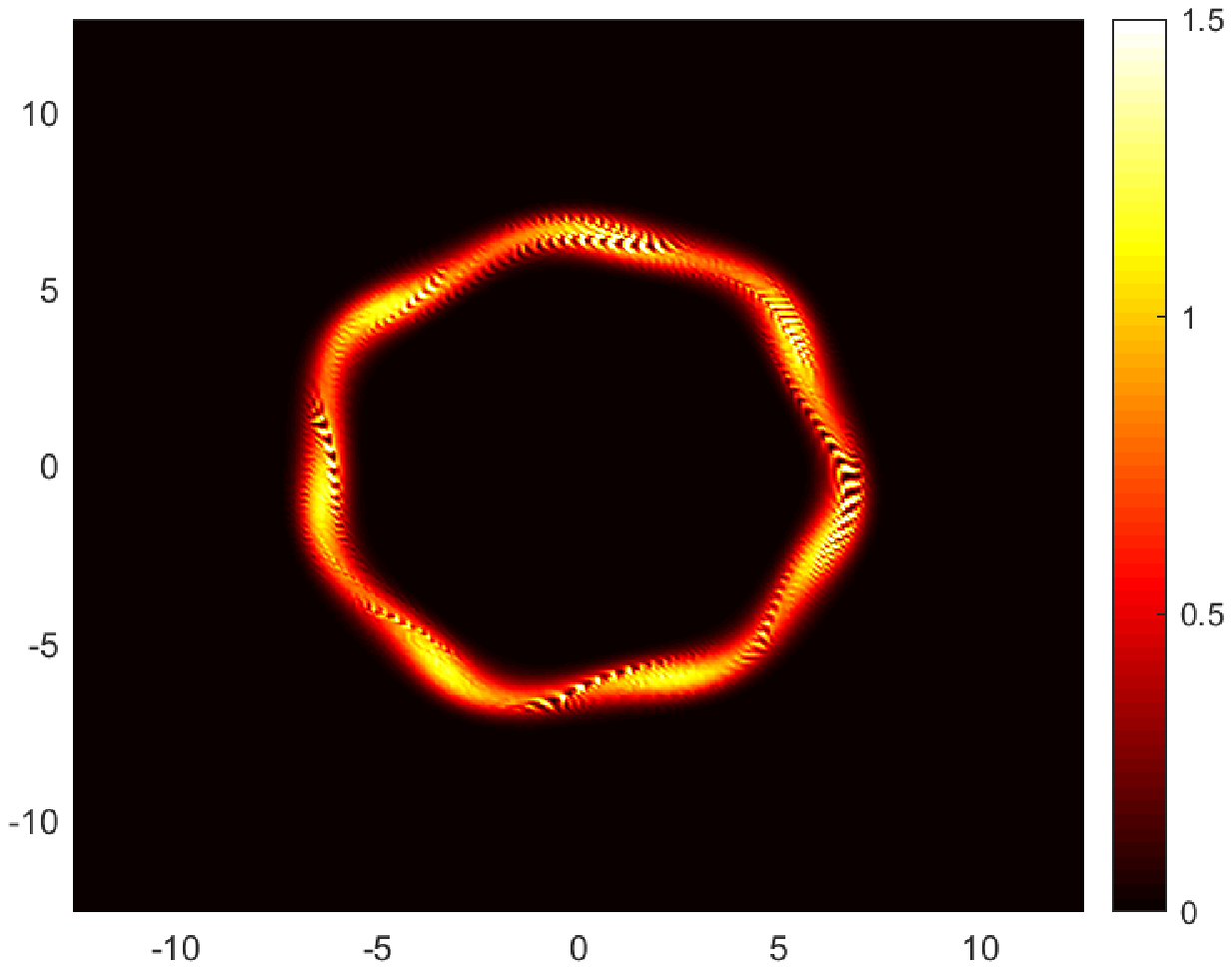}
}
\quad
\subfigure{
\includegraphics[scale=.4]{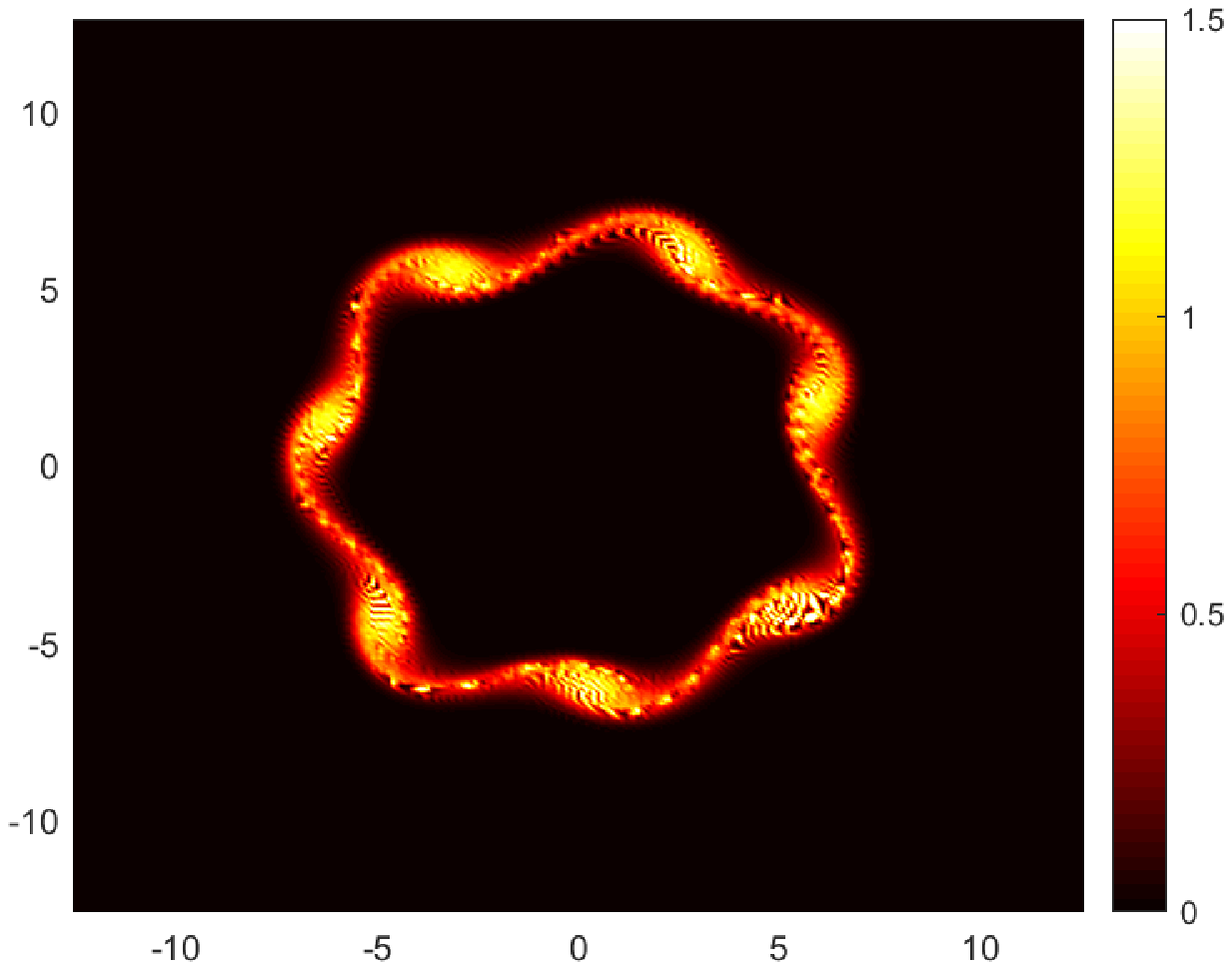}
}
\quad
\subfigure{
\includegraphics[scale=.4]{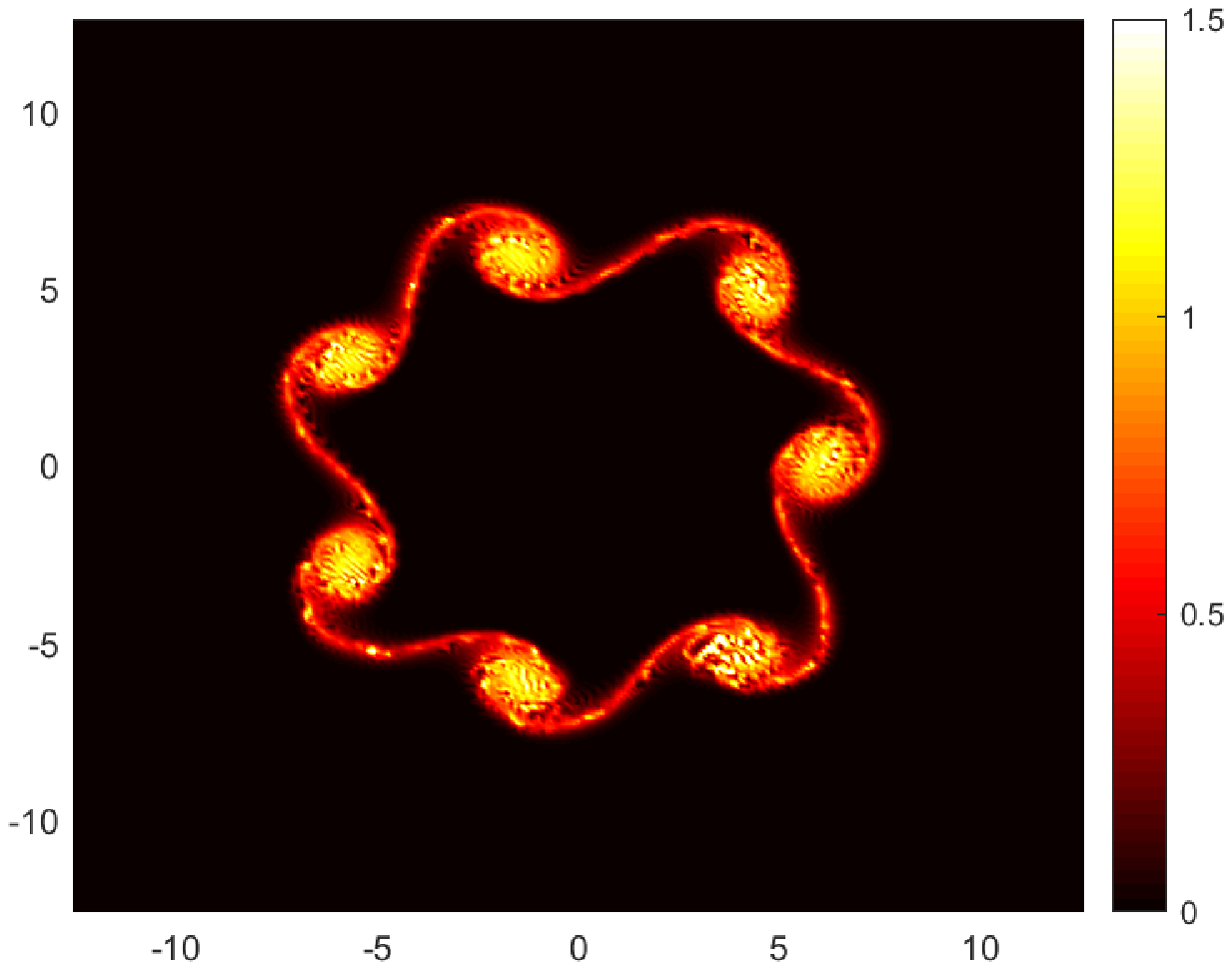}
}
\quad
\subfigure{
\includegraphics[scale=.4]{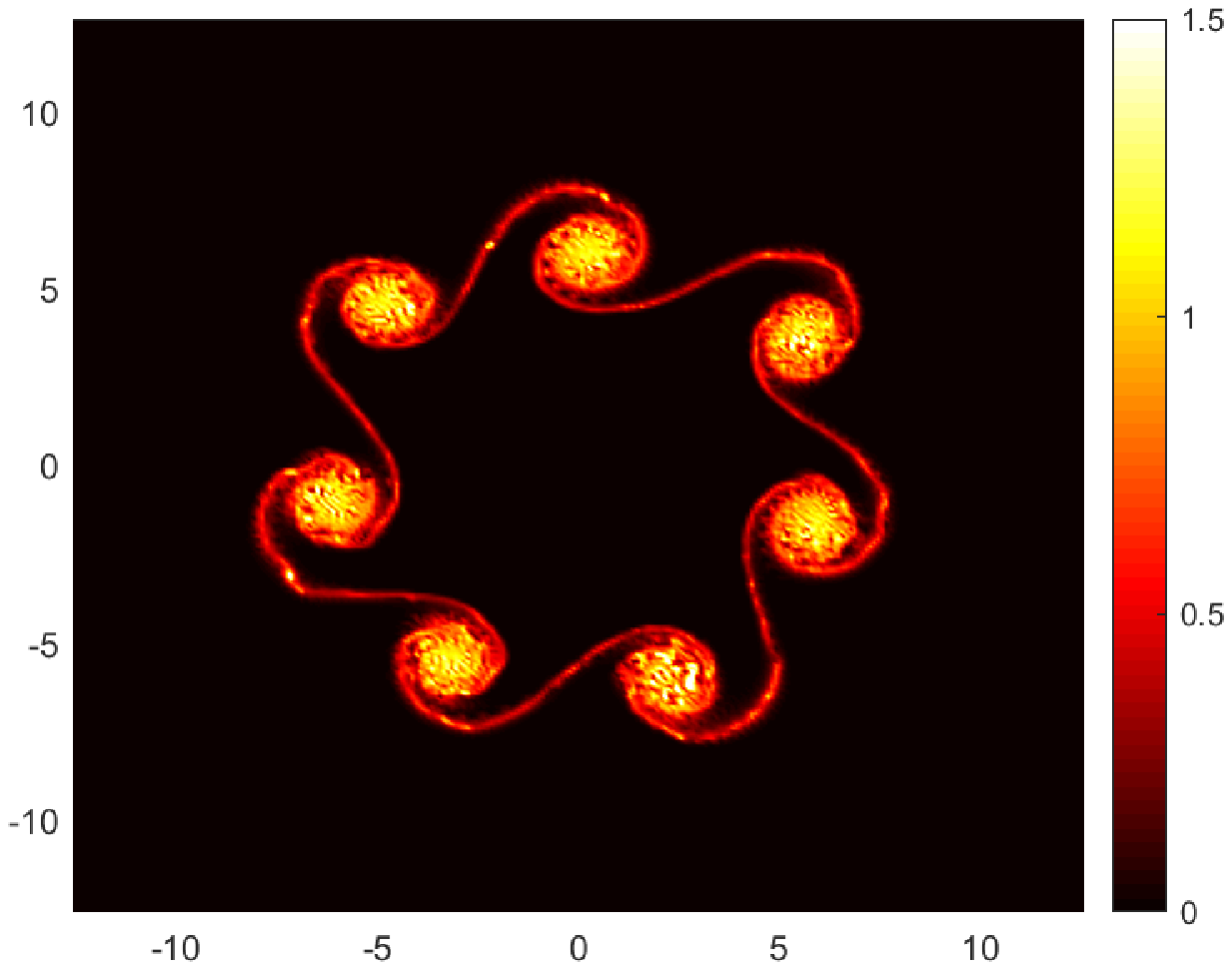}
}
\caption{Time evolution of the density $\rho$ for time $t=5,10,15,20$ and $l=7$ with  strength of magnetic field $\varepsilon=0.01$.}
\label{imgD3}
\end{figure}

Then, we take a small value of  $\varepsilon=0.1$ in Figure \ref{imgD2} and $\varepsilon=0.01$ in Figure \ref{imgD3}. In our computation, we use $\Delta{t}=0.1$ and $\Delta{t}=0.01$ respectively.

We plot the development of the diocotron instability at different time with $l=5$ in Figure \ref{imgD2}. It is observed that the five vortex structures are  well captured by our method. If we change the direction of magnetic field, it can be observed in this figure that the vortex structure will move on the reversed direction.

From Figure \ref{imgD3}, it can be observed that the vortex structures become smaller when the magnetic field is stronger. Also, the image becomes not clear.
To improve this more particles in practical computation are needed and finer resolution should be applied. The computational cost can still be distributed efficiently due to the
parallel approach used in this paper.

\section{Conclusion}

In this work, we have developed the symplectic Particle-in-Cell methods for solving the Vlasov--Poisson system according to its Poisson bracket. We use the appropriate finite element spaces so that the semi-discrete system is equipped with a discrete Poisson structure. With regard to temporal discretization, we apply the splitting technique. As each subsystem can be solved exactly, the resulting discretization can preserve the Poisson structure of the concerning system. In order to implement efficiently our algorithm, parallel computing technique has been designed, and applied to various problems. The parallel efficiency of practical computation has  been verified via the numerical results of this experiments.
We also  show the convergence rates of splitting methods of first and second order which are taken  as a benchmark test. The numerical results perform that  they all match the theoretical order very well.
To verify the advantage of Poisson bracket preserving methods constructed in this paper, we study the preservation of energy both theoretically and numerically.
This guarantees the numerical simulation over long-time.
We also present a 2+2-dimensional  example, in this example  the external magnetic field can be strong.
The error and stability analysis of the derived numerical methods will be reported in the future publications.

\section{Acknowledgments}
This research was supported by the National Natural Science Foundation of China (11771436,11505185), the National Key R\&D Program of China (2017YFE0301704).

\appendix

\section{Jacobi identity of the continuous bracket}

\label{app:1}
In this section we prove the Jacobi identity of the bracket~(\ref{eq:MMWB}) when $\nabla\cdot{\bf B}_0=0$. To simplify the notation, we use $[\cdot,\cdot]$ to replace the continuous bracket $\{\{\cdot,\cdot\}\}$ in this section.

Following the idea in~\citep{Morrison2013A}, we rewrite the continuous bracket~(\ref{eq:MMWB}) as two parts,
\[
[\mathcal{F},\mathcal{G}] =[\mathcal{F},\mathcal{G}]_{xv}+[\mathcal{F},\mathcal{G}]_{B}.
\]
Each part in the form can be written as
\begin{align*}
 &[\mathcal{F},\mathcal{G}]_{xv} (f)=\int_{\Omega_{x}\times\Omega_{v}} f\left\{ \frac{\delta\mathcal{F}}{\delta f},\frac{\delta\mathcal{G}}{\delta f}\right\}_{\mathbf{xv}}d\mathbf{x}d\mathbf{v}\\
 &[\mathcal{F},\mathcal{G}]_{B} (f)=\int_{\Omega_{x}\times\Omega_{v}} f{\bf B}_0\cdot\left(\frac{\partial}{\partial\mathbf{v}}\frac{\delta\mathcal{F}}{\delta f}\times\frac{\partial}{\partial\mathbf{v}}\frac{\delta\mathcal{G}}{\delta f}\right)d\mathbf{x}d\mathbf{v}.
\end{align*}

Then the Jacobi identity reads
\[
\begin{aligned}
[[\mathcal{F},\mathcal{G}] ,\mathcal{H}] +cyc&= \underbrace{[[\mathcal{F},\mathcal{G}]_{xv} ,\mathcal{H}]_{xv}}_{1}+\underbrace{[\mathcal{F},\mathcal{G}]_{xv} ,\mathcal{H}]_{B}}_{2}\\
 & +\underbrace{[[\mathcal{F},\mathcal{G}]_{B} ,\mathcal{H}]_{xv}}_{3}+\underbrace{[[\mathcal{F},\mathcal{G}]_{B} ,\mathcal{H}]_{B}}_{4}\\
 & +cyc,
\end{aligned}
\]
where the symbol $cyc$ means cyclic permutation.

Term 1 vanishes because of the Jacobi identity of $\{\cdot,\cdot\}_{{\bf xv}}$. It is also the bracket of VP equation under $\mathbf{B}_0=\mathbf{0}$, see~\citep{Marsden1999Introduction,Casas2017}.

Next, we simplify the notations. Let $\delta\mathcal{F}=\frac{\delta\mathcal{F}}{\delta f},\delta\mathcal{G}=\frac{\delta\mathcal{G}}{\delta f},\delta\mathcal{H}=\frac{\delta\mathcal{H}}{\delta f}$ and $\alpha=\frac{\partial}{\partial\mathbf{v}}\delta\mathcal{F},\beta=\frac{\partial}{\partial\mathbf{v}}\delta\mathcal{G},\gamma=\frac{\partial}{\partial\mathbf{v}}\delta\mathcal{H}$. And we use ${\bf B}_i$ or $\alpha_i$ to denote the $i$-th component of ${\bf B}_0$ or $\alpha$.
Then term 2+3 reads
\begin{equation*}
\int_{\Omega_{x}\times\Omega_{v}} f{\bf B}_0\cdot(\frac{\partial}{\partial\mathbf{v}}\{\delta\mathcal{F},\delta\mathcal{G}\}_{{\bf xv}}\times\gamma)+f\{{\bf B}_0\cdot(\alpha\times\beta),\delta\mathcal{H}\}_{{\bf xv}}d\mathbf{x}d\mathbf{v}.
\end{equation*}

Notice that $(u\times v)_i=\sum_{jk}\epsilon_{ijk}u_j v_k$ for two vectors $u,v$,
\begin{equation*}
\begin{aligned}
&{\bf B}_0\cdot(\frac{\partial}{\partial\mathbf{v}}\{\delta\mathcal{F},\delta\mathcal{G}\}_{{\bf xv}}\times\gamma)+\{{\bf B}_0\cdot(\alpha\times\beta),\delta\mathcal{H}\}_{{\bf xv}}\\
=&\sum_{ijk}\epsilon_{ijk}{\bf B}_i(\{\alpha_j,\delta\mathcal{G}\}_{{\bf xv}}+\{\delta\mathcal{F},\beta_j\}_{{\bf xv}})\gamma_k+\sum_{ijk}\epsilon_{ijk}\{{\bf B}_i\alpha_j\beta_k,\delta\mathcal{H}\}_{{\bf xv}},
\end{aligned}
\end{equation*}
where $\epsilon_{ijk}$ is the Levi-Civita symbol.
As ${\bf B}_0$ is independent of $\mathbf{v}$, the third term on the right side of the above equality reads
\begin{equation*}
\begin{aligned}
&\{{\bf B}_i\alpha_j\beta_k,\delta\mathcal{H}\}_{{\bf xv}}\\
=&\sum_{l}(\frac{\partial}{\partial\mathbf{x}_l}({\bf B}_i\alpha_j\beta_k)\gamma_l-{\bf B}_i\frac{\partial}{\partial\mathbf{v}_l}(\alpha_j\beta_k)\frac{\partial}{\partial\mathbf{x}_l}\delta\mathcal{H})\\
=&{\bf B}_i\{\alpha_j\beta_k,\delta\mathcal{H}\}_{{\bf xv}}+\sum_{l}\frac{\partial}{\partial\mathbf{x}_l}{\bf B}_i(\alpha_j\beta_k\gamma_l).
\end{aligned}
\end{equation*}

By Leibniz’ rule of $\{\cdot,\cdot\}_{{\bf xv}}$, it is known that
\begin{equation*}
\begin{aligned}
&\{\alpha_j,\delta\mathcal{G}\}_{{\bf xv}}\gamma_k+\{\delta\mathcal{F},\beta_j\}_{{\bf xv}}\gamma_k+\{\alpha_j\beta_k,\delta\mathcal{H}\}_{{\bf xv}}\\
=&\{\alpha_j,\delta\mathcal{G}\}_{{\bf xv}}\gamma_k+\{\delta\mathcal{F},\beta_j\}_{{\bf xv}}\gamma_k+\{\beta_k,\delta\mathcal{H}\}_{{\bf xv}}\alpha_j+\{\alpha_j,\delta\mathcal{H}\}_{{\bf xv}}\beta_k\\
=&\underbrace{(\{\alpha_j,\delta\mathcal{G}\}_{{\bf xv}}\gamma_k+\{\beta_k,\delta\mathcal{H}\}_{{\bf xv}}\alpha_j)}_{\mathcal{A}}+\underbrace{(\{\delta\mathcal{F},\beta_j\}_{{\bf xv}}\gamma_k+\{\alpha_j,\delta\mathcal{H}\}_{{\bf xv}}\beta_k)}_{\mathcal{B}}.
\end{aligned}
\end{equation*}

The terms $\sum_{ijk}\epsilon_{ijk}{\bf B}_i\mathcal{A}$ and $\sum_{ijk}\epsilon_{ijk}{\bf B}_i\mathcal{B}$ cancel out by permuting $\mathcal{F}$, $\mathcal{G}$ and $\mathcal{H}$.
Notice that $\sum_{ijk}\sum_{l}\epsilon_{ijk}\frac{\partial}{\partial\mathbf{x}_l}{\bf B}_i(\alpha_j\beta_k\gamma_l)=\sum_{l}\frac{\partial}{\partial\mathbf{x}_l}{\bf B}_0\cdot(\alpha\times\beta)\gamma_l$ and $(\alpha\times\beta)_i\gamma_j+cyc=\delta_{ij}(\alpha\times\beta)\cdot\gamma$, the summation of terms 2 and 3 leads to $\int_{\Omega_{x}\times\Omega_{v}} f(\nabla\cdot{\bf B}_0)((\alpha\times\beta)\cdot\gamma)d\mathbf{x}d\mathbf{v}$ which vanishes when $\nabla\cdot{\bf B}_0=0$.

According to the above notations, term 4 reads
\begin{equation*}
\int_{\Omega_{x}\times\Omega_{v}} f{\bf B}_0\cdot(\frac{\partial}{\partial\mathbf{v}}({\bf B}_0\cdot(\alpha\times\beta))\times\gamma)d\mathbf{x}d\mathbf{v}
\end{equation*}
By calculation, it is
\begin{equation*}
\begin{aligned}
&{\bf B}_0\cdot(\frac{\partial}{\partial\mathbf{v}}({\bf B}_0\cdot(\alpha\times\beta))\times\gamma)\\
=&\sum_{ijk}\sum_{lmn}\epsilon_{ijk}\epsilon_{lmn}{\bf B}_i{\bf B}_l(\frac{\partial}{\partial\mathbf{v}_j}\alpha_m\beta_n)\gamma_k\\
=&\sum_{ijk}\sum_{lmn}\epsilon_{ijk}\epsilon_{lmn}{\bf B}_i{\bf B}_l(\alpha_{mj}\beta_n\gamma_k+\alpha_{m}\beta_{nj}\gamma_k),
\end{aligned}
\end{equation*}
where $\alpha_{mj}$ means $\frac{\partial}{\partial\mathbf{v}_j}\alpha_m$.

By shifting the indices $i\rightarrow l$, $j\rightarrow m$, $k\rightarrow n$, $l\rightarrow i$, $m\rightarrow k$, $n\rightarrow j$ and $\alpha\beta\gamma\rightarrow \gamma\alpha\beta$, $\epsilon_{ijk}\epsilon_{lmn}{\bf B}_i{\bf B}_l\alpha_{m}\beta_{nj}\gamma_k\rightarrow \epsilon_{lmn}\epsilon_{ikj}{\bf B}_l{\bf B}_i\gamma_{k}\alpha_{jm}\beta_n$.  

Then term 4 vanishes accordingly due to the anti-symmetry of the Levi-Civita symbol.

By the above calculation, we conclude that the continuous bracket~(\ref{eq:MMWB}) is Poisson if $\nabla\cdot \mathbf{B}_{0}=0.$

\section{Derivation of the discrete Poisson bracket}

\label{app:2} 
As follows, we prove that the particle system (\ref{eq:Vlasov-1}) with the approximate electric field (\ref{eq:Eh}) is still a Hamiltonian system with the defined discrete  Poisson bracket.
Under the approximation (\ref{eq:Disf}), we denote a functional of $f_h$ as $$\mathcal{F}[f_h]=F(\mathbf{ \omega,X,V}).$$ Denote
\[
f_s(\mathbf{x},\mathbf{v},t)=\omega_{s}\delta(\mathbf{x}-\mathbf{X}_{s})\delta(\mathbf{v}-\mathbf{V}_{s}), s=1,2,\ldots,N_p.
\]
The discrete variables can be reexpressed by $\omega_{s}=\int f_{s}d\mathbf{x}d\mathbf{v}$,
$\mathbf{X}_{s}=\frac{1}{\omega_{s}}\int\mathbf{x}f_{s}d\mathbf{xd\mathbf{v}}$ and
$\mathbf{V}_{s}=\frac{1}{\omega_{s}}\int\mathbf{v}f_{s}d\mathbf{x}d\mathbf{v}$.
Taking their derivatives w.r.t $f_s$  providing
\begin{gather*}
\frac{\delta\omega_{s}}{\delta f_{s}}=1,\quad\frac{\delta\mathbf{X}_{s}}{\delta f_{s}}=\frac{\mathbf{x}-\mathbf{X}_{s}}{\omega_{s}},\quad\frac{\delta\mathbf{V}_{s}}{\delta f_{s}}=\frac{\mathbf{v}-\mathbf{V}_{s}}{\omega_{s}}.
\end{gather*}
Notice that $\frac{\delta \mathcal{F}}{\delta f_h}=\sum_s\frac{\delta \mathcal{F}}{\delta f_{s}}$.  Using the chain rule of variation  to calculate
$\frac{\delta \mathcal{F}}{\delta f_{s}}$ gives
\begin{equation*}
\begin{aligned}
 \frac{\delta \mathcal{F}}{\delta f_{s}}& = \frac{\delta\omega_{s}}{\delta f_{s}}\frac{\partial F}{\partial\omega_{s}}+\frac{\delta\mathbf{X}_{s}}{\delta f_{s}}\frac{\partial F}{\partial\mathbf{X}_{s}}+\frac{\delta\mathbf{V}_{s}}{\delta f_{s}}\frac{\partial F}{\partial\mathbf{V}_{s}}\\
 & = \frac{\partial F}{\partial\omega_{s}}+\frac{\mathbf{x}-\mathbf{X}_{s}}{\omega_{s}}\frac{\partial F}{\partial\mathbf{X}_{s}}+\frac{\mathbf{v-}\mathbf{V}_{s}}{\omega_{s}}\frac{\partial F}{\partial\mathbf{V}_{s}}.
\end{aligned}
\label{eq:variationalfs}
\end{equation*}
It leads to
\begin{align*}
\left\{ \frac{\delta\mathcal{F}}{\delta f_{s}},\frac{\delta\mathcal{G}}{\delta f_{s}}\right\} _{{\bf xv}} & =\left\{ \frac{\partial F}{\partial\omega_{s}}+\frac{\mathbf{x}-\mathbf{X}_{s}}{\omega_{s}}\frac{\partial F}{\partial\mathbf{X}_{s}}+\frac{\mathbf{v-}\mathbf{V}_{s}}{\omega_{s}}\frac{\partial F}{\partial\mathbf{V}_{s}},\frac{\partial G}{\partial\omega_{s}}+\frac{\mathbf{x}-\mathbf{X}_{s}}{\omega_{s}}\frac{\partial G}{\partial\mathbf{X}_{s}}+\frac{\mathbf{v-}\mathbf{V}_{s}}{\omega_{s}}\frac{\partial G}{\partial\mathbf{V}_{s}}\right\} _{{\bf xv}}\\
 & =\frac{1}{\omega_{s}^{2}}\left(\frac{\partial F}{\partial\mathbf{X}_{s}}\cdot\frac{\partial G}{\partial\mathbf{V}_{s}}-\frac{\partial F}{\partial\mathbf{V}_{s}}\cdot\frac{\partial G}{\partial\mathbf{X}_{s}}\right),
\end{align*}
Then
\[
\begin{aligned}
\int_{\Omega_x\times\Omega_v}  f_h\left\{ \frac{\delta\mathcal{F}}{\delta f_h},\frac{\delta\mathcal{G}}{\delta f_h}\right\}_{\mathbf{xv}}d\mathbf{x}d\mathbf{v}
=\sum_{s=1}^{N_p}\frac{1}{\omega_{s}}\left(\frac{\partial F}{\partial\mathbf{X}_{s}}\cdot\frac{\partial G}{\partial\mathbf{V}_{s}}-\frac{\partial F}{\partial\mathbf{V}_{s}}\cdot\frac{\partial G}{\partial\mathbf{X}_{s}}\right).\end{aligned}
\]
And
\[
\begin{aligned}
\int_{\Omega_x\times\Omega_v} f_h{\bf B}_0(\mathbf{x})\cdot\left(\frac{\partial}{\partial\mathbf{v}}\frac{\delta\mathcal{F}}{\delta f_h}\times\frac{\partial}{\partial\mathbf{v}}\frac{\delta\mathcal{G}}{\delta f_h}\right)d\mathbf{x}d\mathbf{v}
& =\int_{\Omega_x\times\Omega_v} f_h{\bf B}_0(\mathbf{x})\cdot\left(\frac{1}{\omega_{s}^{2}}\frac{\partial F}{\partial\mathbf{V}_{s}}\times\frac{\partial G}{\partial\mathbf{V}_{s}}\right)d\mathbf{x}d\mathbf{v}\\
 & =\sum_{s=1}^{N_p} \frac{1}{\omega_{s}}{\bf B}_0({\bf X}_s)\cdot\left(\frac{\partial F}{\partial\mathbf{V}_s} \times\frac{\partial G}{\partial\mathbf{V}_s}\right).
\end{aligned}
\]

Therefore, we can define the following discrete bracket operator as
\begin{equation*}
\begin{aligned}
\left\{ F,G\right\} (\mathbf{X},\mathbf{V})
=\sum_{s=1}^{N_p}\frac{1}{\omega_{s}}\left(\frac{\partial F}{\partial\mathbf{X}_{s}}\cdot\frac{\partial G}{\partial\mathbf{V}_{s}}-\frac{\partial G}{\partial\mathbf{X}_{s}}\cdot\frac{\partial F}{\partial\mathbf{V}_{s}}\right)+\sum_{s=1}^{N_p} \frac{1}{\omega_{s}}{\bf B}_0(\mathbf{X}_s)\cdot\left(\frac{\partial F}{\partial\mathbf{V}_s} \times\frac{\partial G}{\partial\mathbf{V}_s}\right).
\end{aligned}
\end{equation*}

\section{Jacobi identity of the discrete bracket}

\label{app:3}
As follows, we prove that the discrete bracket (\ref{eq:dPoisson}) is a Poisson bracket under $\nabla\cdot\mathbf{B}_{0}(\mathbf{x})=0$.

At the end of section 4, we connect the discrete Poisson bracket (\ref{eq:dPoisson}) which derivation has been shown in \ref{app:2} to a matrix $\mathbb{K}$ which is called the Poisson matrix. Thus, the discrete bracket is Poisson if and only if for the component of matrix $\mathbb{K}$ the following identity holds for all $i,j,k$
\begin{equation*}
\sum_{l=1}^{6{N_p}}(\frac{\partial\mathbb{K}_{ij}}{\partial\mathbf{Z}^l}\mathbb{K}_{lk}+\frac{\partial\mathbb{K}_{jk}}{\partial\mathbf{Z}^l}\mathbb{K}_{li}+\frac{\partial\mathbb{K}_{ki}}{\partial\mathbf{Z}^l}\mathbb{K}_{lj})=0,
\end{equation*}
where $\mathbf{Z}^l$ is the $l$-th component of $\mathbf{Z}$, see~\citep{Olver1986,Hairer2006}.

It can be known that 
\begin{align*}
&\sum_{l=1}^{6{N_p}}(\frac{\partial\mathbb{K}_{ij}}{\partial\mathbf{Z}^l}\mathbb{K}_{lk}+\frac{\partial\mathbb{K}_{jk}}{\partial\mathbf{Z}^l}\mathbb{K}_{li}+\frac{\partial\mathbb{K}_{ki}}{\partial\mathbf{Z}^l}\mathbb{K}_{lj})\\
=&\sum_{l=1}^{3{N_p}}(\frac{\partial\mathbb{N}^{-1}\hat{\mathbb{B}}({\bf X})_{ij}}{\partial\mathbf{X}^l}\mathbb{N}^{-1}_{lk}+\frac{\partial\mathbb{N}^{-1}\hat{\mathbb{B}}({\bf X})_{jk}}{\partial\mathbf{X}^l}\mathbb{N}^{-1}_{li}+\frac{\partial\mathbb{N}^{-1}\hat{\mathbb{B}}({\bf X})_{ki}}{\partial\mathbf{X}^l}\mathbb{N}^{-1}_{lj}).
\end{align*}
Due to that  $\mathbb{N}$ is a diagonal constant matrix, the right term of the  above equality leads to
\begin{align*}&\frac{\partial\mathbb{N}^{-1}\hat{\mathbb{B}}({\bf X})_{ij}}{\partial\mathbf{X}^k}\mathbb{N}^{-1}_{kk}+\frac{\partial\mathbb{N}^{-1}\hat{\mathbb{B}}({\bf X})_{jk}}{\partial\mathbf{X}^i}\mathbb{N}^{-1}_{ii}+\frac{\partial\mathbb{N}^{-1}\hat{\mathbb{B}}({\bf X})_{ki}}{\partial\mathbf{X}^j}\mathbb{N}^{-1}_{jj}.
\end{align*}

If $\mathbf{X}^i$ and $\mathbf{X}^j$ are not the components for the same particle position $\mathbf{X}_s$, then $\frac{\partial\mathbb{N}^{-1}\hat{\mathbb{B}}({\bf X})_{ij}}{\partial\mathbf{X}^k}\mathbb{N}^{-1}_{kk}+\frac{\partial\mathbb{N}^{-1}\hat{\mathbb{B}}({\bf X})_{jk}}{\partial\mathbf{X}^i}\mathbb{N}^{-1}_{ii}+\frac{\partial\mathbb{N}^{-1}\hat{\mathbb{B}}({\bf X})_{ki}}{\partial\mathbf{X}^j}\mathbb{N}^{-1}_{jj}=0$  due to that  the value of all partial derivatives of functions is zero. Therefore, we only need to consider the case in which  $\mathbf{X}^i,\mathbf{X}^j,\mathbf{X}^k$ are the components for the same particle position $\mathbf{X}_s$.
This reads
\begin{align*}
\frac{1}{\omega^2_s}(\frac{\partial\hat{\mathbb{B}}_{ij}}{\partial\mathbf{X}^k}+\frac{\partial\hat{\mathbb{B}}_{jk}}{\partial\mathbf{X}^i}+\frac{\partial\hat{\mathbb{B}}_{ki}}{\partial\mathbf{X}^j})=0.
\end{align*}
In the above equality, if the two of three indexes are equal the terms related  cancel out  because of the skew-symmetry of matrix $\hat{\mathbb{B}}$. With the Levi-Civita symbol $\epsilon_{ijk}$, it is known that $\hat{\mathbb{B}}_{ij}$ has the expression $\hat{\mathbb{B}}_{ij}=\epsilon_{ijk}{\bf B}_0(\mathbf{X}_s)_{k}$ when the three indices are different. Then $\frac{\partial\hat{\mathbb{B}}_{ij}}{\partial\mathbf{X}^k}+\frac{\partial\hat{\mathbb{B}}_{jk}}{\partial\mathbf{X}^i}+\frac{\partial\hat{\mathbb{B}}_{ki}}{\partial\mathbf{X}^j}=0$ if and only if $\nabla_{\mathbf{X}_{s}}\cdot \mathbf{B}_{0}(\mathbf{X}_{s})=0$. This implies the bracket defined by $\mathbb{K}$ is Poisson if $\nabla\cdot \mathbf{B}_{0}(\mathbf{x})=0.$

\bibliographystyle{plain}
\bibliography{refs}

\end{document}